\documentclass[12pt]{article}
\usepackage{amsmath,amssymb}

\newtheorem{theorem}{Theorem}[section]
\newtheorem{proposition}{Proposition}[section]
\newtheorem{corollary}{Corollary}[section]  
\newtheorem{lemma}{Lemma}[section]
  
\newcommand{\CQFD}{\nolinebreak\hfill\rule{2mm}{2mm}\medbreak\par}  
\newtheorem{Rem}{Remark}[section]  
  
\numberwithin{equation}{section}

\def\bbe{\mathbb E}

\def\bbp{\mathbb P}
\def\bbr{\mathbb R}
\def\udel{\Delta}

\DeclareMathOperator{\textvar}{Var}
\DeclareMathOperator{\textcov}{Cov}

\begin{document}

\title{On the Longest Increasing Subsequence for Finite and Countable 
Alphabets}
\author{Christian Houdr\'e \thanks{Georgia Institute of Technology, 
School of Mathematics, Atlanta, Georgia, 30332-0160
{\it E-mail address:} houdre@math.gatech.edu} 
\and Trevis J. Litherland \thanks{Georgia Institute of Technology, 
School of Mathematics, Atlanta, Georgia, 30332-0160
{\it E-mail address:} trevisl@math.gatech.edu} }

\maketitle

\vspace{0.5cm}

\begin{abstract}
\noindent Let $X_1, X_2, \ldots, X_n, \ldots $ be a sequence of  
iid random variables with values in a finite alphabet $\{1,\ldots,m\}$.
Let $LI_n$ be the length of the longest increasing
subsequence of $X_1, X_2, \ldots, X_n.$  
We express the limiting distribution of $LI_n$ as
functionals of $m$ and $(m-1)$-dimensional Brownian motions. 
These expressions are then related to similar functionals
appearing in queueing theory, allowing us to further 
establish asymptotic behaviors as $m$ grows.  The finite 
alphabet results are then used to treat the countable (infinite) 
alphabet.    

\end{abstract}

\noindent{\footnotesize {\it AMS 2000 Subject Classification:} 60C05, 60F05, 60F17, 
60G15, 60G17, 05A16}
%05A16: Asymptotic enmeration
%60C05: Combinatorial probability
%60F05: Central limit and other weak theorems
%60F17: Functional limit theorems
%60G15: Gaussian processes
%60G17: Sample path properties

\noindent{\footnotesize {\it Keywords:} Longest increasing subsequence, 
Brownian functional, Functional Central Limit Theorem, 
Tracy-Widom distribution.}

\section{Introduction}

    The pursuit of a robust understanding of the 
asymptotics of the length of the longest increasing
subsequence $L\sigma_n$ of a random permutation of length $n$ 
-- often known as "Ulam's Problem" --
has given rise to a remarkable collection of results 
in recent years.   
The work of
Logan and Shepp \cite{LS}, and Vershik and Kerov \cite{VK}, 
first showed that $ {\bbe L\sigma_n}/{\sqrt{n}} \rightarrow 2$.
Following this fundamental asymptotic result,
Baik, Deift, and Johansson,
in their landmark paper \cite{BDJ}, determined
the limiting distribution of $L\sigma_n$, 
properly centered and normalized. 
This problem has emerged as 
a nexus of once seemingly unconnected mathematical ideas. 
Indeed, the latter paper is, in particular, 
quite remarkable for the sheer breadth of 
mathematical machinery required, machinery calling upon an
understanding of random matrix theory, the asymptotics of Toeplitz operators, 
Riemann-Hilbert Theory, as well as the Robinson-Schensted-Knuth correspondence, 
to obtain the limiting Tracy-Widom distribution.  

Initial approaches to the problem relied heavily
on combinatorial arguments.  
Most work of the last decade, however, such as
that of Aldous and Diaconis \cite{AD1}
and Sepp\"al\"ainen \cite{S2}, have
instead used interacting particle
processes and so-called "hydrodynamical arguments"
to show that $L\sigma_n/\sqrt{n} \rightarrow 2$ 
in expectation and in probability. 
Building on these ideas, Groeneboom \cite{Gr2}
proves such convergence results 
using only the convergence of random signed measures,
while Cator and Groeneboom \cite{CG1} prove that
$ {\bbe L\sigma_n}/{\sqrt{n}} \rightarrow 2$ in
a way that avoids both ergodic decomposition arguments
and the subadditive ergodic theorem.
Aldous and Diaconis \cite{AD2} also connect these 
particle process concepts to the card game patience sorting.  
Finally, Sepp\"al\"ainen \cite{S} employs these
particle processes to a verify an open asymptotics
problem in queueing theory.
Moving beyond the asymptotics of $\bbe L\sigma_n$,
Cator and Groeneboom \cite {CG2} use particle processes to directly
obtain the cube-root asymptotics of the variance
of $L\sigma_n$.  Further non-asymptotic results for
$L\sigma_n$ are found in \cite{Gr}.

The related problem of the asymptotics of $LI_n$ when 
the sequence is drawn uniformly from a finite
alphabet of size $m$ has developed along parallel lines.  
Tracy and Widom \cite{TW}, as well as Johannson \cite{Jo} ,
have shown that the limiting distribution 
again enjoys a direct connection to the 
distribution of the largest eigenvalue 
in the Gaussian Unitary Ensemble paradigm.  
Its, Tracy, and Widom \cite{ITW1,ITW2} have
further examined this problem in the inhomogeneous 
case, relating the limiting distribution
to certain direct sums of GUEs.  In another direction, 
Chistyakov and G\"otze \cite{ChG}
have pursued the two-letter Markov case.

Problems from statistical physics have long 
inspired a lot of the research into these topics.
Kuperberg \cite{Ku}, for instance, shows that certain quantum spin
matrices are, in law, asymptotically equal to a traceless GUE matrix.
The standard general overview of the 
subject of random matrices is Mehta \cite{Me},
a work motivated and influenced by 
some of the origins of the subject in physics.

While the above achievements have undoubtedly stimulated further inquiry, 
one might still suspect that a more direct
route to the limiting distribution of 
$LI_n$ might be had, one whose methods reflect the essentially
probabilistic nature of the problem.  This paper 
proposes a step towards such an approach 
for the independent finite alphabet case, 
calling only upon some very well-known results of classical probablity theory
described below.  Indeed, the sequel will 
show that the limiting distribution of $LI_n$ 
can be constructed in a most natural manner
as a Brownian functional.  In the context of
random growth processes, Gravner, Tracy, and Widom \cite{GTW}
have already obtained a Brownian functional of the form
we derive.  This functional appeared first in the
work of Glynn and Whitt \cite{GW}, in queueing
theory, and its relation to the eigenvalues of the GUE
has independently been studied by
Baryshnikov \cite{Ba}.
It is, moreover, remarked in \cite {GTW} 
that the longest increasing 
subsequence problem could also
be studied using a Brownian functional formulation.

We begin our study of this problem, in 
the next section, by expressing $LI_n$ as
a simple algebraic expression.  Using 
this simple characterization, we then briefly
determine, in Section $3$, 
the limiting distribution 
of $LI_n$ (properly centered and normalized) in the case 
of an $m$-letter alphabet 
with each letter drawn 
independently.  Our result is expressed
as a functional of an 
$(m-1)$-dimensional Brownian motion
with correlated coordinates.
Using certain natural symmetries, 
this limiting distribution
is further expressed as various 
functionals of a (standard) Brownian motion.  
In Section $4$, connections with
the Brownian functional originating 
with the work of Glynn and Whitt in
queueing theory are investigated.
This allows us to investigate 
the asymptotics as $m$ grows.  Section $5$ is devoted to 
obtaining the corresponding results for countable alphabets.  
In Section $6$, we finish 
the paper by indicating
some open questions and future 
directions for research.

\section{Combinatorics}

Let $X_1, X_2, \ldots, X_n, \ldots $  consist of a sequence of values taken from an
$m$-letter ordered alphabet, $\alpha_1 < \alpha_2 < \cdots < \alpha_m$.  
Let $a^r_k$ be the number of occurrences
of $\alpha_r\in \{1,\dots,m\}$ among $X_1,X_2,\dots, X_k$, $1 \le k \le n$.  
Each increasing subsequence of $X_1,X_2,\dots, X_n$ consists simply of runs
of identical values, with the values of each successive run forming 
an increasing subsequence of $\alpha_r$.
Moreover, the number of occurrences of $\alpha_r\in \{\alpha_1,
\dots,\alpha_m\}$ among
$X_{k+1},\dots, X_\ell$, where $1 \le k < \ell \le n$, is simply
$a^r_{\ell}-a^r_k$. The length of the longest increasing subsequence 
of $X_1,X_2,\dots, X_n$ is then given by

\begin{equation}\label{item1}
 LI_n=\max_{\stackrel{\scriptstyle 0\le k_1\le\cdots}{\le k_{m-1}\le n}}
[(a^1_{k_1}-a^1_0)+(a^2_{k_2}-a^2_{k_1})+\cdots +
(a^m_n-a^m_{k_{m-1}})],\end{equation}
{\it i.e.},
\begin{equation}\label{item2}
 LI_n=\max_{\stackrel{\scriptstyle 0\le k_1\le\cdots}{\le k_{m-1}\le n}}
[(a^1_{k_1}-a^2_{k_1})+(a^2_{k_2}-a^3_{k_2})+\cdots +
(a^{m-1}_{k_{m-1}}-a^m_{k_{m-1}})+a^m_n],\end{equation}
where $a^r_0=0$.  
For $i = 1, \dots ,n$ and $r = 1, \dots ,{m-1}$, let 
\begin{equation}\label{item3}
Z^r_i=\begin{cases} 1, &\text{if $X_i=\alpha_r,$}\\
-1, & \text{if $X_i=\alpha_{r+1},$}\\
0, & \text{otherwise,}\end{cases}
\end{equation}
and let $S^r_k=\sum^k_{i=1}Z^r_i$, $k = 1, \dots ,n$, with also $S^r_0=0$. 
Then clearly $S^r_k = a^r_k - a^{r+1}_k$.  Hence,
\begin{equation}\label{item4}
LI_n=\max_{\stackrel{\scriptstyle 0\le k_1\le\cdots}{\le k_{m-1}\le n}}
\{S^1_{k_1}+S^2_{k_2}+ \cdots + S^{m-1}_{k_{m-1}}+a^m_n\}.\end{equation}

Since $a^1_k, \dots ,a^m_k$ must evidently sum to $k$, we have

\begin{align*}
n &= \sum^m_{r=1}a^r_n \\
&= -\sum^{m-1}_{r=1}r(a^{r+1}_n-a^r_n) + ma^m_n\\
&= -\sum^{m-1}_{r=1}rS^r_n + ma^m_n.\\
\end{align*}

Solving for $a^m_n$ gives us
$$a^m_n=\frac nm- \frac1m \sum^{m-1}_{r=1} rS^r_n.$$

Substituting into \eqref{item4}, we finally obtain
\begin{equation}\label{item5}
 LI_n=\frac nm-\frac1m\sum^{m-1}_{r=1} rS^r_n+
\max_{\stackrel{\scriptstyle 0\le k_1\le\cdots}{\le k_{m-1}\le n}}
\{S^1_{k_1}+S^2_{k_2}+ \cdots + S^{m-1}_{k_{m-1}}\}.\end{equation}

The expression \eqref{item5} is of a {\it purely
combinatorial nature or, in more probabilistic terms, is of a pathwise nature}.  
We now analyze \eqref{item5} in light
of the probabilistic nature of the sequence $X_1,X_2,\dots, X_n$.

\section{Probabilistic Development}

Throughout the sequel, Brownian functionals will play a central r\^ole.  
By a {\it Brownian motion} we shall mean an a.s.~continuous, centered 
Gaussian process $B(t)$, $0 \le t \le 1$, with $B(0) = 0$, 
having stationary, independent increments.  
By a {\it standard Brownian motion} we shall
mean that Var$B(t) = t$, $0 \le t \le 1$, {\it i.e.}, 
we endow  $C[0,1]$ with the Wiener measure. 
A {\it standard $m$-dimensional Brownian 
motion} will be defined to be 
a vector-valued process consisting of  
$m$ independent Brownian motions.
More generally, an {\it $m$-dimensional Brownian motion} 
shall refer to a linear transformation of
a standard $m$-dimensional Brownian motion.
Throughout the paper, we assume that our underlying
probability space is rich enough so that all 
the Brownian motions and sequences we study
can be defined on it.

We consider first the case in which 
$X_1, X_2, \ldots, X_n, \ldots $  are iid, with each
letter drawn uniformly from 
${\cal A} = \{\alpha_1,\dots, \alpha_m\}.$  
Then for each fixed letter $r$, the sequence 
$Z^r_1,Z^r_2,\dots\,Z^r_n,\dots$ is also formed 
of iid random variables with $\bbp(Z^r_1=1) = \bbp(Z^r_1=-1) = 1/m$, 
and $\bbp(Z^r_1=0) = 1 - 2/m$.

Thus $\bbe Z^r_1=0$, and $\bbe(Z^r_1)^2=2/m$, and so, 
Var$S^r_n=2n/m$, for $r = 1,2, \dots, {m-1}.$  
Defining $\hat B^r_n(t)=\frac1{\sqrt{2n/m}}
S^r_{[nt]}+\frac1{\sqrt{2n/m}} (nt-[nt])Z^r_{[nt]+1}$, 
for $0 \le t \le 1$, and noting that the local
maxima of $\hat B^i_n(t)$ occur at 
$t=k/n$, $k=0,\dots,n$, we have from \eqref{item5} that

\begin{equation}\label{item6}
 \frac{LI_n-n/m}{\sqrt{2n/m}} = 
-\frac1m\sum^{m-1}_{i=1} i\hat B^i_n(1)+
\max_{\stackrel{\scriptstyle 0\le t_1\le\cdots}{\le t_{m-1}\le 1}}
[\hat B^1_n(t_1)+\cdots + \hat B^{m-1}_n (t_{m-1})].\end{equation}

We note that we can now invoke Donsker's Theorem 
since the measures $\bbp_n$
generated by $(\hat B^1_n(t),\dots, \hat B^{m-1}_n(t))$ 
satisfy $\bbp_n(A) \rightarrow \bbp_{\infty}(A)$, for all 
Borel subsets $A$ of the space of 
continuous functions  $C([0,1]^{m-1})$ for which
$\bbp_{\infty}(\partial A) = 0$, where $\bbp_{\infty}$
is the limiting $(m-1)$-dimensional Wiener measure.
Then, applying Donsker's Theorem and 
the Continuous Mapping Theorem we have that 
$(\hat B^1_n(t),\dots, \hat B^{m-1}_n(t))\Rightarrow (\tilde{B}^1(t),\dots,
\tilde{B}^{m-1}(t))$, where the Brownian motion 
on the right has a covariance structure
which we now describe.  First, Cov$(Z^r_1,Z^s_1) = \bbe Z^r_1 Z^s_1=0$, 
for $|r-s|\ge 2$, and Cov$(Z^r_1,Z^{r+1}_1) =\bbe Z^r_1 Z^{r+1}_1=-1/m$,
for $r = 1,2, \dots, {m-1}$.  Then, as already noted, 
for each fixed $r$, $Z^r_1,Z^r_2,\dots\, Z^r_n,\dots$ 
are iid, and for fixed $k$, 
$Z^1_k, Z^2_k, \dots, Z^{m-1}_k$ are 
dependent but identically distributed random
variables.  Moreover, it is equally 
clear that for any $r$ and $s$, $1 \le r < s \le m-1$, 
the sequences $(Z^r_k)_{k \ge 1}$ and $(Z^s_{\ell})_{{\ell} \ge 1}$
are also identical distributions of the $Z^r_k$
and that $Z^r_k$ and $Z^s_{\ell}$ are independent for $k \neq {\ell}$.
Thus, Cov$(S^r_n,S^s_n) = n$Cov$(Z^r_1,Z^s_1)$.
This result, together with our $2n/m$ 
normalization factor gives the following 
covariance matrix for $(\tilde{B}^1(t),\dots,\tilde{B}^{m-1}(t))$:

\begin{equation}\label{mattri}
\quad
t\begin{pmatrix} 1 & -1/2 &&&\bigcirc\\
-1/2 & 1 & -1/2\\
&\ddots &\ddots & \ddots\\
\bigcirc &&\ -1/2 & 1 & -1/2\\
&&& -1/2 & 1\end{pmatrix}.\\
\end{equation}

We remark here that the functional in 
\eqref{item6} is a bounded linear functional
on $C(0,1)^{m-1}$.  
(This fact will be used throughout the paper.)
Hence, by a final application of the Continuous Mapping Theorem,

\begin{equation}\label{item6b}
\frac{LI_n-n/m}{\sqrt{2n/m}} \Rightarrow 
-\frac1m\sum^{m-1}_{i=1} i \tilde{B}^i(1)+
\max_{\stackrel{\scriptstyle 0\le t_1\le\cdots}{\le t_{m-1}\le 1}}
\sum^{m-1}_{i=1} \tilde{B}^i(t_i).
\end{equation}

We have thus obtained the limiting distribution 
of $LI_n$ as a Brownian functional.
Tracy and Widom \cite{TW} already obtained the 
limiting distribution of $LI_n$ in terms of 
the distribution of the largest eigenvalue of 
the Gaussian Unitary Ensemble (GUE) of $m \times m$ 
Hermitian matrices having trace zero.  
Johansson \cite{Jo} generalized this work to encompass all
$m$ eigenvalues.  Gravner, Tracy, and Widom
\cite{GTW} in their study of random growth processes
make a connection between the distribution of the largest
eigenvalue in the $m \times m$ GUE and a Brownian
functional essentially equivalent,
up to a normal random variable, to the right hand side
of \eqref{item6b}.
(This will become clear as we refine our
understanding of \eqref{item6b} in the sequel.)
For completeness, we now state our result.

\begin{proposition}\label{prop0}
Let $X_1, X_2, \ldots, X_n, \ldots $ be a sequence of iid random
variables drawn uniformly from the ordered finite alphabet 
${\cal A} = \{\alpha_1,\dots,\alpha_m\}$.  
Then
\begin{equation}\label{item6x}
\frac{LI_n-n/m}{\sqrt{2n/m}} \Rightarrow 
-\frac1m\sum^{m-1}_{i=1} i \tilde{B}^i(1)+
\max_{\stackrel{\scriptstyle 0\le t_1\le\cdots}{\le t_{m-1}\le 1}}
\sum^{m-1}_{i=1} \tilde{B}^i(t_i),
\end{equation}

\noindent where $(\tilde{B}^1(t),\dots,\tilde{B}^{m-1}(t))$ is 
an $(m-1)$-dimensional Brownian
motion with covariance matrix given by \eqref{mattri}.
\end{proposition}

For $m=2$, \eqref{item6x} simply becomes

\begin{equation}\label{item6xa}
\frac{LI_n-n/2}{\sqrt{n}} \Rightarrow 
-\frac{1}{2}B(1)+ \max_{{\scriptstyle 0\le t \le 1}} B(t),
\end{equation}

\noindent where $B$ is standard one-dimensional 
Brownian motion.
A well-known result of Pitman \cite{P} 
implies that, up to a factor of $2$,
the functional in \eqref{item6xa} is identical 
in law to the radial part of a three-dimensional 
standard Brownian motion at time $t=1$.  
Specifically, Pitman shows that the process 
$2\max_{0 \le s \le t}B(s) - B(t)$  
is identical in law to 
$\sqrt{(B^1(t))^2+(B^2(t))^2+(B^3(t))^2}$,
where $(B^1(t),B^2(t),B^3(t))$ is 
a standard $3$-dimensional Brownian motion.

Let us now show that the functional in \eqref{item6xa} 
does indeed have the same distribution as
that of the largest eigenvalue of a 
$2 \times 2$ zero-trace matrix of the form

$$ \quad
\begin{pmatrix} X & Y + iZ\\
Y - iZ & -X\end{pmatrix},$$\\

\noindent where $X$, $Y$, and $Z$ are centered 
independent normal random variables, all with variance $1/4$.  
These random variables have a joint density given by

$$ f_3(x,y,z) = \left(\frac{2}{\pi}\right)^{3/2} e^{-2x^2 - 2y^2 - 2z^2}, 
\quad  \text{$(x,y,z) \in \bbr^3$}.$$

It is straightforward to show that the largest 
eigenvalue of our matrix is given by $\lambda_1 =
\sqrt{X^2 + Y^2 + Z^2}.$  Thus,
up to a scaling factor of $2$, $\lambda_1$
is equal in law to the radial Brownian motion 
expression of Pitman at $t=1$.  
Explicitly, since $4\lambda_1^2 =
4X^2 + 4Y^2 + 4Z^2$ consists of the 
sum of the squares of three iid standard normal random variables, 
$4\lambda_1^2$ must have a $\chi^2$ 
distribution with $3$ degrees of freedom.  Since this distribution
has a density of $h(x) = (1/\sqrt{2\pi})x^{1/2}e^{-x/2}$, 
we immediately find that $\lambda_1$ has density 

\begin{align}
g(\lambda_1) &= \frac{1}{\sqrt{2\pi}}(4\lambda_1^2)^{1/2}
e^{-(4\lambda_1^2)/2} (8\lambda_1)\nonumber\\
&= \frac{16}{\sqrt{2\pi}}\lambda_1^2 e^{-2\lambda_1^2}, \quad \text{$\lambda_1 > 0$}.\nonumber
\end{align}

Let us look now at the connection between the $2 \times 2$ GUE
and the traceless matrix we have just analyzed. 
Consider the $2 \times 2$ matrix

$$ \quad
\begin{pmatrix} X_1 & Y + iZ\\
Y - iZ & X_2\end{pmatrix},$$\\

\noindent where $X_1$, $X_2$, $Y$, and $Z$ are 
independent normal random variables, with Var$X_1 = $
Var$X_2 = 1/2$, and with Var$Y = $ Var$Z =1/4$.  
Since these random variables have a joint
density given by

$$ f_4(x_1,x_2,y,z) = \frac{2}{\pi^2}e^{-x_1^2 - x_2^2 - 2y^2 - 2z^2}, 
\quad  \text{$(x_1,x_2,y,z) \in \bbr^4$},$$

\noindent conditioning on the zero-trace subspace 
$\left\{X_1 + X_2 = 0\right\}$, and using the transformation
$X_1^\prime = (X_1 - X_2)/\sqrt{2}$ and
$X_2^\prime = (X_1 + X_2)/\sqrt{2}$,
we obtain the conditional density

\begin{align*}
f_3(x_1^\prime,y,z) = &\left(\frac{2}{\pi}\right)^{3/2} e^{-2(x_1^\prime)^2 - 2y^2 - 2z^2}
\end{align*}

\noindent which is also the joint density of three iid centered 
normal random variables $X_1^\prime$, $Y$, and $Z$ with common variance $1/4$,
which we had previously obtained.

Let us finally note that one can directly evaluate \eqref{item6xa} 
in a classical manner using the Reflection Principle
to obtain the corresponding density (see, {\it e.g.} \cite{GTW,HLM}). \\

It is instructive to express \eqref{item6x} 
in terms of an $(m-1)$-dimensional
standard Brownian motion $(B^1(t),\dots, B^{m-1}(t))$.  
It is not hard to check that we can express
$\tilde{B}^i(t)$, $i=1,\dots,{m-1}$, 
in terms of the $B^i(t)$ as follows:

\begin{equation}\label{item6y}
\tilde{B}^i(t) =
\begin{cases}
B^1(t),  & \text{$i=1$},\\
\sqrt{\frac{i+1}{2i}}B^i(t)-\sqrt{\frac{i-1}{2i}}B^{i-1}(t), & \text{$2\le i\le m-1$}.

\end{cases}
\end{equation}

Substituting \eqref{item6y} back into \eqref{item6x}, we obtain a more
symmetric expression for our limiting distribution:

\begin{equation}\label{item6z}
\frac{LI_n-n/m}{\sqrt{n}} \Rightarrow 
\frac{1}{\sqrt{m}}  \max_{\stackrel{\scriptstyle 0\le t_1\le\cdots}
{\le t_{m-1}\le t_m = 1}}
\sum^{m-1}_{i=1} \!\! \left[-\sqrt{\frac{i}{\!i+1}}B^i(t_{i+1}) 
+ \sqrt{\frac{i+1}{i}}B^i(t_i)\right]\!.
\end{equation}

The above Brownian functional is similar to one 
introduced by Glynn and Whitt \cite{GW}, 
in the context of a queueing problem:

\begin{equation}\label{item6za}
D_m = 
\max_{\stackrel{\scriptstyle 0 = t_0 \le t_1\le\cdots}
{\le t_{m-1}\le t_m  = 1}}
\sum^m_{i=1} \left[ B^i(t_i) - B^i(t_{i-1}) \right],
\end{equation}

\noindent where $(B^1(t),\dots, B^m(t))$ is 
an $m$-dimensional standard Brownian motion.
Gravner, Tracy, and Widom \cite{GTW}, in 
studying a one-dimensional discrete space and 
discrete time process, have shown that its 
limiting distribution is equal in law to both
$D_m$ and also to the largest eigenvalue 
$\lambda_1^{(m)}$ of an $m \times m$  
Hermitian matrix taken from a GUE.  
That is, $D_m$  and $\lambda_1^{(m)}$ are 
in fact identical in law.  Independently,
Baryshnikov \cite{Ba}, studying closely 
related problems of queueing theory and of 
monotonous paths on the integer lattice, 
has shown that the {\it process}
$(D_m)_{m \ge 1}$ has the same law as 
the {\it process} $(\lambda_1^{(m)})_{m \ge 1}$, 
where $\lambda_1^{(m)}$ is the largest
eigenvalue of the matrix consisting of
the first $m$ rows and $m$ columns of
an infinite matrix in the Gaussian Unitary Ensemble.

\begin{Rem}It is quite clear that $LI_n \ge n/m$ 
a.s., since at least one of the $m$ letters must
lie on a substring of length at least $n/m$.  
Hence, the limiting functional in \eqref{item6x}
must be supported on the positive real line.  
Can we see directly that a.s.~the functional on the right
hand side of \eqref{item6z} is also non-negative?  
Indeed, for consider the more general Brownian
functional of the form 

$$\max_{\stackrel{\scriptstyle 0\le t_1\le\cdots}{\le t_{m-1}\le t_m = 1}}
\sum^{m-1}_{i=1} \left[ \beta_i B^i(t_{i+1}) - \eta_i B^i(t_i) \right],$$

\noindent where $0 \le \beta_i \le \eta_i,$ 
for $i = 1,2,\dots,m-1$.  Now for any
fixed $ t_{i+1} \in (0,1]$, $i=1,\dots,{m-1}$,
$\max_{0 \le t_i \le t_{i+1}} 
\left[ \beta_i B^i(t_{i+1}) - \eta_i B^i(t_i) \right]$
is at least as large as the maximum 
value at the two extremes, that is, 
when $t_i=0$ or $t_i = t_{i+1}$.
These two values are simply 
$\beta_i B^i(t_{i+1})$ and $(\beta_i-\eta_i) B^i(t_{i+1})$.  
Since $0 \le \beta_i \le \eta_i$ a.s., at 
least one of these two values is 
non-negative.  Hence, we can successively
find $t_{m-1}, t_{m-2},\dots,t_1$ such 
that {\it each} term of the functional is non-negative a.s.
Thus the whole functional must be non-negative a.s.  Taking
$\beta_i = \sqrt{i/(i+1)}$ and $\eta_i = \sqrt{(i+1)/i}$, 
the result holds for \eqref{item6z}.
The functional of Glynn and Whitt in 
\eqref{item6za} does not succumb to the same 
analysis since the $i=1$ term demands that $t_0=0$.
\end{Rem}

Let us now turn our attention to the 
$m$-letter case wherein each
letter $\alpha_r\in $ occurs with 
probability $0<p_r<1$, independently,
and the $p_r$ need not be equal as in the previous uniform case.
For the non-uniform case, 
Its, Tracy, and Widom in \cite{ITW1} and \cite{ITW2} 
obtained the limiting distribution of $LI_n$.
Reordering the probabilities such that
$p_1 \ge p_2 \ge \cdots \ge p_m$,
and grouping those probabilities
having identical values
$p_{(j)}$ of multiplicity $k_j$, $j = 1,\dots,d$, 
(so that $\sum_{j=1}^d k_j = m$
and $\sum_{j=1}^d p_{(j)}k_j = 1$), 
they show that the limiting 
distribution is identical in law 
to the distribution
of the largest eigenvalue of 
the direct sum of $d$ mutually independent 
$k_j \times k_j$ GUEs, whose eigenvalues
$(\lambda_1,\lambda_2,\dots,\lambda_m) = 
(\lambda_1^{k_1}, \lambda_2^{k_1}, \dots, \lambda_{k_1}^{k_1},\dots,
\lambda_1^{k_d}, \lambda_2^{k_d}, \dots, \lambda_{k_d}^{k_d})$
satisfy $\sum_{i=1}^m\sqrt{p_i}\lambda_i = 0$.
With the above ordering of the probabilities,
the limiting distribution simplifies to a $k_1$-fold
integral involving only $p_1$ and $k_1$. 
(See Remark~\ref{sumrmt} for some explicit
expressions and more details.)
We now state our own result in terms of functionals of
Brownian motion.

\begin{theorem}\label{thm2}
Let $X_1, X_2, \ldots, X_n, \ldots $ be a sequence of iid random
variables such that $\bbp(X_1 = \alpha_r) = p_r,$ for  $r = 1, \ldots, m$, where
$0 < p_r < 1$ and $\sum^m_{r=1}p_r = 1.$ Then

\begin{equation}\label{item6zz}
\frac{LI_n-p_{max}n}{\sqrt n} \Rightarrow - \frac1m
\sum^{m-1}_{i=1} i\sigma_i \tilde{B}^i(1)+
\mathop{\max_{\stackrel{\scriptstyle 0 = t_0 \le t_1\le \cdots}
{\le t_{m-1}\le t_m = 1}}}_{
t_i=t_{i-1},\ i\in I^*}
\sum^{m-1}_{i=1} \sigma_i \tilde{B}^i(t_i),
\end{equation}

\noindent where $p_{max}=\max_{1\le r\le m}p_r$,
$\sigma^2_r = p_r+p_{r+1}-(p_r-p_{r+1})^2,$ 
$I^*=\{r \in \{1,\ldots,m\}:p_r<p_{max}\}$,
and where $(\tilde{B}^1(t),\dots,\tilde{B}^{m-1}(t))$ 
is an $(m-1)$-dimensional Brownian
motion with covariance matrix given by
$$ t\begin{pmatrix}
1 & \rho_{1,2} & \rho_{1,3} & \cdots & \rho_{1,m-1}\\
\rho_{2,1} &1  & \rho_{2,3} & \cdots & \rho_{2,m-1}\\
\vdots & \vdots & \ddots & \ddots & \vdots\\
\vdots & \vdots && 1 &\rho_{m-2,m-1} \\
\rho_{m-1,1} &\rho_{m-1,2} &\cdots &\rho_{m-1,m-2} &1\end{pmatrix},$$

\noindent with 
$$\rho_{r,s}=\begin{cases}
-p_r-\frac{\mu_r\mu_s}{\sigma_r\sigma_{s}}, & \text{$s=r-1,$}\\
-p_s-\frac{\mu_r\mu_s}{\sigma_r\sigma_{s}}, & \text{$s=r+1,$}\\
- \frac{\mu_r\mu_s}{\sigma_r\sigma_s}, & \text{$|r-s| > 1$, 
\quad$1\le r, s\le m-1$,}\end{cases}$$

\noindent and with $\mu_r = p_r - p_{r+1}, 1 \le r \le m-1.$

\end{theorem}

\noindent \begin{Proof} As before, we 
begin with the expression for $LI_n$ 
displayed in \eqref{item5},
noting that for each letter $\alpha_r$,  
$1 \le r \le m-1$, $(Z^r_k)_{k \ge 1}$ 
forms a sequence of iid random
variables, and that moreover $Z^r_k$ 
and $Z^s_{\ell}$ are independent for $k \neq {\ell}$,
and for any $r$ and $s$. Now, however, for 
each fixed $k$, the $Z^r_k$ are no longer identically 
distributed; indeed,

\begin{equation}\label{item7}
\begin{cases}
\mu_r:=\bbe Z^r_1=p_r-p_{r+1}, & \text{$1\le r\le m-1$,}\\
\sigma^2_r:= \mbox{Var} Z^r_1 = p_r+p_{r+1}-(p_r-p_{r+1})^2, & \text{$1\le r\le m-1$}.
\end{cases}
\end{equation}

Since $0 < p_r < 1$, we have 
$\sigma^2_r > 0$ for all $1 \le r \le m-1$.
We are thus led to define our 
approximation to a Brownian motion by
$$\hat B^r_n(t):= \frac{S^r_{[nt]}-\mu_r n}{\sigma_r\sqrt n} + (nt-[nt])\frac{Z^r_{[nt]+1}}{\sigma_r\sqrt n}, 
\qquad 0 \le t \le 1,  \quad r=1, \dots, {m-1}. $$  

Again noting that the local maxima of $\hat B^i_n(t)$ occur 
on the set $\{t:t=k/n, k=0,\dots,n\}$, 
\eqref{item5} becomes

\begin{align}\label{item8}
LI_n=\frac nm&-\frac1m\sum^{m-1}_{i=1} i\left[\sigma_i\hat B^i_n(1)
\sqrt n+\mu_in\right]\nonumber\\
&+
\max_{\stackrel{\scriptstyle 0 = t_0 \le t_1\le \cdots}{\le t_{m-1}\le t_m = 1}}
\left\{\sum^{m-1}_{i=1}\left[\sigma_i\hat B^i_n(t_i)\sqrt n+
\mu_it_in\right]\right\}.\end{align}

Next,
\begin{align*}
\sum^{m-1}_{i=1}i\mu_i&=
\sum^{n-1}_{i=1}\sum^{m-1}_{j=i}\mu_j=
\sum^{m-1}_{i=1}\sum^{m-1}_{j=i} (p_j-p_{j+1})\\
&=\sum^{m-1}_{i=1} (p_i-p_m)= (1-p_m)-(m-1)p_m\\
&=1-mp_m.
\end{align*}

Hence, \eqref{item8} becomes
\begin{align}\label{item9}
LI_n=\frac nm&- \frac{(1-m p_m)n}m - \frac1m\sum^{m-1}_{i=1} i\sigma_i\hat B^i_n(1)
\sqrt n \nonumber\\
&+
\max_{\stackrel{\scriptstyle 0 = t_0 \le t_1\le \cdots}{\le t_{m-1}\le t_m = 1}}
\sum^{m-1}_{i=1}\left[\sigma_i\hat B^i_n(t_i)\sqrt n+
\mu_it_in\right], \end{align}
and, dividing through by $\sqrt n$, we obtain
\begin{align}\label{item10}
\frac{LI_n}{\sqrt n} &= p_m\sqrt n- \frac1m
\sum^{m-1}_{i=1} i\sigma_i\hat B^i_n(1)\nonumber\\
&\qquad +
\max_{\stackrel{\scriptstyle 0 = t_0 \le t_1\le \cdots}{\le t_{m-1}\le t_m = 1}}
\sum^{m-1}_{i=1} \left[\sigma_i\hat B^i_n(t_i)+\mu_it_i\sqrt n\right].
\end{align}

Let $t_0=0$, and let 
$\udel _i=t_i-t_{i-1}$, $i=1,\dots, m-1$. Since

$$\sum^{m-1}_{i=1}\mu_it_i
=\sum^{m-1}_{i=1}\mu_i\sum^{i}_{j=1}\udel_i
=\sum^{m-1}_{i=1}\udel_i\sum^{m-1}_{j=i}\mu_j
=\sum^{m-1}_{i=1}\udel _i(p_i-p_m),$$

\noindent \eqref{item10} becomes

\begin{align}\label{item11}
\frac{LI_n}{\sqrt n} = p_m\sqrt n&- \frac1m
\sum^{m-1}_{i=1} i\sigma_i\hat B^i_n(1)\nonumber\\
&+
\max_{\stackrel{\scriptstyle \udel_i\ge 0}{\sum^{m-1}_{i=1}\udel_i\le 1}}
\left\{\sum^{m-1}_{i=1} \sigma_i\hat B^i_n(t_i)+\sqrt n
\sum^{m-1}_{i=1}\udel_i(p_i-p_m)\right\},
\end{align}

\noindent where $t_i=\sum^i_{j=1}\udel _j$.

Setting also $\udel_m=1-t_{m-1}$ ({\it i.e.}, $t_m:=1$), 
\eqref{item11} enjoys a more symmetric
representation as

\begin{align}\label{item12}
\frac{LI_n}{\sqrt n} = &- \frac1m
\sum^{m-1}_{i=1} i\sigma_i\hat B^i_n(1)\nonumber\\
&+
\max_{\stackrel{\scriptstyle \udel_i\ge 0}{\sum^{m-1}_{i=1}\udel_i= 1}}
\left[\sum^{m-1}_{i=1} \sigma_i\hat B^i_n(t_i)+\sqrt n
\sum^{m}_{i=1}\udel_ip_i\right].
\end{align}

Let $p_{max}=\max_{1\le i\le m}p_i$.  Then
\begin{align}\label{item13}
\frac{LI_n-p_{max}n}{\sqrt n} = &- \frac1m
\sum^{m-1}_{i=1} i\sigma_i\hat B^i_n(1)\nonumber\\
&+
\max_{\stackrel{\scriptstyle \udel_i\ge 0}{\sum^{m-1}_{i=1}\udel_i= 1}}
\left[\sum^{m-1}_{i=1} \sigma_i\hat B^i_n(t_i)+\sqrt n
\sum^{m}_{i=1}\udel_i(p_i-p_{max})\right].
\end{align}

Clearly, if $\udel_i>0$ for any $i$ such that $p_i<p_{max}$, then
$$\sqrt n\sum^m_{i=1}\udel_i (p_i-p_{max})\mathop{\longrightarrow}^{a.s.} - \infty.$$

Intuitively, then, we should demand that
$\udel_i=0$ for $i\in I^*:=\{i \in \{1,2,\dots,m\}:p_i<p_{max}\}$.  Indeed,
we now show that in fact

\begin{equation}\label{item14}
\frac{LI_n-p_{max}n}{\sqrt n} = - \frac1m
\sum^{m-1}_{i=1} i\sigma_i\hat B^i_n(1)+
\mathop{\max_{\stackrel{\scriptstyle 0 = t_0 \le t_1\le \cdots}
{\le t_{m-1}\le t_m = 1}}}_{
t_i=t_{i-1},\ i\in I^*}
\sum^{m-1}_{i=1} \sigma_i\hat B^i_n(t_i) + E_n,
\end{equation}

\noindent where the remainder term $E_n$ is a random variable 
converging to zero in probability 
as $n \rightarrow \infty $.

To see this, let us introduce the following notation.  Writing \\
$t = (t_1,t_2,\dots,t_{m-1})$,
let $T = \{t: 0 \le t_1 \le \cdots \le t_{m-1} \le 1\}$ and let
$T^* = \{t \in T: t_i = t_{i-1}, i \in I^*\}$.  Setting 
$C_n(t) = \sum^{m-1}_{i=1} \sigma_i\hat B^i_n(t_i)$ 
and $R(t) = \sum^{m}_{i=1}(t_i-t_{i-1})(p_{max}-p_i)$,
we can, respectively, rewrite the maximal 
terms of \eqref{item13} and \eqref{item14} 
as $$\max_{t \in T} \left[C_n(t) - \sqrt n R(t)\right]$$ and
$$\max_{t \in T^*} C_n(t).$$ 

\noindent By the compactness of 
$T$ and $T^*$ and the continuity of
$C_n(t)$ and $R(t)$, we see that for 
each $n$ and each $\omega \in \Omega$, there is a
$\tau_n \in T$ and a $\tau^*_n \in T^*$ such that 
$$C_n(\tau_n) - \sqrt n R(\tau_n) 
= \max_{t \in T} \left[C_n(t) - \sqrt n R(t)\right],$$ and
$$C_n(\tau^*_n) = \max_{t \in T^*} C_n(t).$$ 
(Note that the piecewise-linear nature 
of $C_n(t)$ and the linear nature of $R(t)$ imply
that the arguments maximizing the above must 
lie on a finite set and that the measurablility 
of $\tau_n$ and $\tau^*_n$ is trivial.)

Now we first claim that the set of 
optimizing arguments $\{\tau_n\}^{\infty}_{n=1}$ 
a.s.~does not have an accumulation point 
lying outside of $T^*$.  Suppose the contrary,
namely that for each $\omega$ 
in a set $A$ of positive measure,
there is a subsequence 
$(\tau_{n_k})^{\infty}_{k=1}$ of 
$(\tau_n)^{\infty}_{n=1}$ such 
that $d(\tau_{n_k},T^*) > \epsilon$,
for some $\epsilon > 0$, where the 
metric $d$ is the one induced by the
$L_{\infty}$-norm over $T$, {\it i.e.}, by
$\|t\|_{\infty} = \max_{1\le i\le m-1}|t_i|$. 

Then, since $T^* \subset T$, 
it follows that, for all $n$,
$$C_{n}(\tau_{n}) - \sqrt{n} R(\tau_{n}) \ge C_{n}(\tau^*_{n}),$$
almost surely. Now if $p_{max} = p_m$, 
then $t = (0,\dots,0) \in T^*$, and
if for some $1 \le j \le m-1$ we have 
$p_{max} = p_j > \max_{j+1 \le i \le m} p_i$,   
then $t = (0,\dots,0,1,\dots,1) \in T^*$, 
where there are $j$
zeros in $t$. Hence 
$C_{n_k}(\tau^*_{n_k}) \ge  C_{n_k}(0,\dots,0,1,\dots,1)
=\sum^{m-1}_{i=j+1} \sigma_i\hat B^i_{n_k}(1)$ a.s.,
where the sum is taken to be zero for $j = m$.   
Given $0 < \delta < 1$, by the 
Central Limit Theorem, we can find a
sufficiently negative real $\alpha$ such that

\begin{align*}
\bbp \left(C_{n_k}(\tau_{n_k}) - \sqrt{n_k} R(\tau_{n_k}) 
\ge \alpha\right) &\ge \bbp \left(C_{n_k}(\tau^*_{n_k}) \ge \alpha \right) \\
&\ge \bbp \left(\sum^{m-1}_{i=j+1} \sigma_i\hat B^i_{n_k}(1) \ge \alpha  \right) \\
&> 1 - \delta,
\end{align*}

\noindent for $n_k$ large enough.  
In particular, this implies that

\begin{equation}\label{item14b}
\bbp \left(A \cap \{C_{n_k}(\tau_{n_k}) - \sqrt{n_k} R(\tau_{n_k})  
\ge \alpha\} \right) > \frac12\bbp(A),
\end{equation}

\noindent for $n_k$ large enough.

Next, note that for any $t \in T$, 
we can modify its components $t_i$ 
to obtain an element of $T^*$, by collapsing 
certain consecutive $t_i$s to single values,
where $i \in \{j-1,j,\dots,\ell\}$ 
and $\{j,j+1,\dots,\ell\} \subset I^*$.  
With this observation, it is not 
hard to see that by replacing such
maximal consecutive sets of 
components $\{t_i\}^{\ell}_{i=j-1}$ 
with their median values, we must have 

$$d(\tau_{n_k},T^*) = \max_{ \{(j,\ell):\{j,j+1,\dots,\ell\} 
\subset I^* \} } \frac{(\tau^\ell_{n_k} - \tau^{j-1}_{n_k})}{2}.$$

Writing $p_{(2)}$ for the largest of the $p_i < p_{max}$, 
we see that for all $k$, and for almost all $\omega \in A$,

\begin{align*}
R(\tau_{n_k}) &= \sum^{m}_{i=1}(\tau^i_{n_k}-\tau^{i-1}_{n_k})(p_{max}-p_i)\\
&=   \sum_{i\in I^*}(\tau^i_{n_k}-\tau^{i-1}_{n_k})(p_{max}-p_i)\\
&\ge (p_{max}-p_{(2)})\sum_{i\in i^*}(\tau^i_{n_k}-\tau^{i-1}_{n_k})\\
&\ge 2 (p_{max}-p_{(2)})d(\tau_{n_k},T^*) \ge 2(p_{max}-p_{(2)})\epsilon.
\end{align*}\\

Now by Donsker's Theorem 
and the Continuous Mapping Theorem, we have that 
$$\max_{t\in T^*}\sum^{m-1}_{i=1} \sigma_i\hat B^i_{n_k}(t_i) 
\Rightarrow \max_{t\in T^*}\sum^{m-1}_{i=1} \sigma_i\tilde{B}^i_{n_k}(t_i),$$
as $n_k \rightarrow \infty$, where 
$(\tilde{B}^1(t),\dots,\tilde{B}^{m-1}(t))$ 
is an $(m-1)$-dimensional Brownian motion
described in greater detail below.  
The point here is simply that 
this limiting functional exists.  
Hence, given $0 < \delta < 1$, 
if $M$ is chosen large enough, then

\begin{align*}
\bbp \left(C_{n_k}(\tau_{n_k}) \le M\right) &
\ge \bbp \left(C_{n_k}(\tau^*_{n_k}) \le M\right) \\
&= \bbp \left(\max_{t\in T^*}\sum^{m-1}_{i=1} 
\sigma_i\hat B^i_{n_k}(t_i) \le M  \right)\\
&> 1 - \delta,
&\end{align*}

\noindent for $n_k$ large enough.  

We now can see how the boundedness 
of $R(\tau_{n_k})$ on $A$ influences that of the
whole expression 
$C_{n_k}(\tau_{n_k}) - \sqrt{n_k} R(\tau_{n_k})$ 
by the following estimates.  
Given $M > 0$ as above, if $k$ is large enough, then
$$n_k \ge \sqrt{(M-\alpha+1)/(2(p_{max}-p_{(2)})\epsilon)},$$ 
and also

\begin{align*}
&\bbp \left(A \cap \{C_{n_k}(\tau_{n_k}) - \sqrt{n_k} R(\tau_{n_k}) 
\le \alpha-1\} \right)\\
  &\qquad = \bbp \left(A \cap \{ C_{n_k}(\tau_{n_k}) \le \alpha-1  + 
\sqrt{n_k} R(\tau_{n_k})\}\right)\\
  &\qquad \ge  \bbp \left(A \cap \{ C_{n_k}(\tau_{n_k}) \le \alpha -1  + 
\sqrt{n_k} (2(p_{max}-p_{(2)})\epsilon)\}\right)\\
  &\qquad \ge  \bbp \left(A \cap \{ C_{n_k}(\tau_{n_k}) \le M\} \right)\\
  &\qquad >  \frac12\bbp(A).
\end{align*}

\noindent But this contradicts \eqref{item14b}, and our optimal
parameter sequences $(\tau_n)^{\infty}_{n=1}$ 
must a.s.~have their accumulation points in $T^*$.

Thus, given $\epsilon > 0$, 
there is an integer $N_\epsilon$ such that the set 
$A_{n,\epsilon} = \{d(\tau_k,T^*) < \epsilon^3, k \ge n\}$ satisfies 
$\bbp(A_{n,\epsilon}) \ge 1-\epsilon$,
for all $n \ge N_\epsilon$. Now for each 
$\tau_n$ define $\hat\tau_n \in T^*$ to be the 
(not necessarily unique) point of $T^*$ which is 
closest in the $L^{\infty}$-distance to $\tau_n$. 
Recalling that 
$$E_n = C_n(\tau_n) - \sqrt n R(\tau_n) - C_n(\tau^*_n) \ge 0,$$ 
almost surely, and noting that 
$R(t) \ge 0$, for all $t \in T$, 
we can estimate the remainder term 
$E_n $ as follows: for $n \ge N_\epsilon$,

{\allowdisplaybreaks
\begin{align}\label{item14c}
\bbp \left(E_n \ge \epsilon  \right) &=  \bbp \left(\{E_n \ge \epsilon\} 
\cap  A_{n,\epsilon} \right) + \bbp \left(\{E_n \ge \epsilon\} 
\cap  A^c_{n,\epsilon} \right)\nonumber\\
&\le \bbp \left(\{E_n \ge \epsilon\} \cap  A_{n,\epsilon} \right) + 
\bbp \left( A^c_{n,\epsilon} \right)\nonumber\\
&<   \bbp \left(\{E_n \ge \epsilon\} \cap  A_{n,\epsilon} \right) + 
\epsilon \nonumber\\
&=   \bbp \left(\{C_n(\tau_n) - \sqrt n R(\tau_n) - C_n(\tau^*_n)   
\ge \epsilon\} \cap  A_{n,\epsilon} \right) + \epsilon \nonumber\\
&\le \bbp \left(\{C_n(\tau_n) - \sqrt n R(\tau_n) - C_n(\hat\tau_n) 
\ge \epsilon\} \cap  A_{n,\epsilon} \right) + \epsilon \nonumber\\
&\le \bbp \left(\{C_n(\tau_n) - C_n(\hat\tau_n) \ge \epsilon\} \cap  A_{n,\epsilon} \right) 
+ \epsilon \nonumber\\
&\le  \bbp \left(  \left|\sum^{m-1}_{i=1} \sigma_i (\hat B^i_n(\tau^i_n) - 
\hat B^i_n(\hat\tau^i_n)) \right| \ge \epsilon \right) + \epsilon. 
\end{align}}

To further bound \eqref{item14c}, note that 
for all $n \ge 1$ and all $1 \le i \le m-1$, we have 
$\mbox{Var}(\hat B^i_n(t_i) - \hat B^i_n(s_i)) = |t_i-s_i|$.
Then, let $(s,t) \in T \times T$ 
be such that $\|t-s\|_{\infty} \le \epsilon^3$.
Using the Bienaym\'e-Chebyshev inequality, 
we find that for $n$ large enough,

\begin{align*}
\bbp \left(  \left|\sum^{m-1}_{i=1} \sigma_i (\hat B^i_n(t_i) - 
\hat B^i_n(s_i)) \right| \ge \epsilon \right)
  &\le \epsilon^{-2}(m-1)^2 \max_{1\le i \le m-1} \sigma_i^2 \|t-s\|_{\infty} \\
  &\le \epsilon^{-2}(m-1)^2 \max_{1\le i \le m-1} \sigma_i^2 \epsilon^3\\
  &=   \epsilon (m-1)^2 \max_{1\le i \le m-1} \sigma_i^2.
\end{align*}

Since $\|\tau_n - \hat \tau_n\| < \epsilon^3$, 
for  $n \ge N_\epsilon$, this can be used to
bound \eqref{item14c}:

\begin{align*}
\bbp \left(\left|E_n\right| \ge \epsilon  \right) &< 
  \bbp \left(  \left|\sum^{m-1}_{i=1} \sigma_i (\hat B^i_n(\tau^i_n) - 
\hat B^i_n(\hat\tau^i_n)) \right| \ge \epsilon \right) + \epsilon\\
& \le  \epsilon \left\{(m-1)^2 \max_{1\le i \le m-1} \sigma_i^2 + 1\right\}.
\end{align*}

\noindent  Finally, $\epsilon$ being arbitrary, 
we have indeed shown that 
$E_n \rightarrow 0$ in probability.\\

Applying Donsker's Theorem, the Continuous Mapping Theorem, 
and the converging together lemma to \eqref{item14}
we finally have:
\begin{equation}\label{item15}
\frac{LI_n-p_{max}n}{\sqrt n} \Rightarrow - \frac1m
\sum^{m-1}_{i=1} i\sigma_i \tilde{B}^i(1)+
\mathop{\max_{\stackrel{\scriptstyle 0 = t_0 \le t_1\le \cdots}
{\le t_{m-1}\le t_m = 1}}}_{
t_i=t_{i-1},\ i\in I^*}
\sum^{m-1}_{i=1} \sigma_i \tilde{B}^i(t_i),
\end{equation}

\noindent where $(\tilde{B}^1(t),\dots,\tilde{B}^{m-1}(t))$ is an 
$(m-1)$-dimensional Brownian motion
 with the following covariance matrix:
$$ t\begin{pmatrix}
1 & \rho_{1,2} & \rho_{1,3} & \cdots & \rho_{1,m-1}\\
\rho_{2,1} & 1 & \rho_{2,3} & \cdots & \rho_{2,m-1}\\
\vdots & \vdots &\ddots & \ddots & \vdots\\
\vdots & \vdots && 1 &\rho_{m-2,m-1} \\
\rho_{m-1,1} &\cdots &\cdots &\rho_{m-1,m-2} &1\end{pmatrix},$$
where 
$$\rho_{r,s}=\begin{cases}
-p_r-\frac{\mu_r\mu_s}{\sigma_r\sigma_{s}}, & \text{$s=r-1,$}\\
-p_s-\frac{\mu_r\mu_s}{\sigma_r\sigma_{s}}, & \text{$s=r+1,$}\\
- \frac{\mu_r\mu_s}{\sigma_r\sigma_s}, & \text{$|r-s| > 1$, 
\quad$1\le r, s\le m-1$.}\end{cases}$$

Now for $t=\ell /n$, and $1 \le r \le s \le m-1$,
the covariance structure above is computed
as follows:

{\allowdisplaybreaks
\begin{align*}
\mbox{Cov}(\hat B^r_n(t),\hat B^s_n(t))&=
\mbox{Cov}\left(\sum^\ell_{i=1}\frac{Z^r_i-\mu_r}{\sigma_r \sqrt n}\,,
\sum^\ell_{i=1}\frac{Z^s_i-\mu_s}{\sigma_s \sqrt n}\right)\\
&=\frac1{n\sigma_r\sigma_s}\mbox{Cov} \left(\sum^\ell_{i=1} 
(Z^r_i-\mu_r),\sum^\ell_{i=1}(Z^s_i-\mu_s)\right)\\
&=\frac1{n\sigma_r\sigma_s}\sum^\ell_{i=1}\mbox{Cov} (Z^r_i-\mu_j,
Z^s_i-\mu_k)\\
&=\frac\ell{n\sigma_r\sigma_s}\mbox{Cov}(Z^r_1-\mu_r, Z^s_1-\mu_s)\quad\\
&=t\begin{cases}\frac 1{\sigma_r\sigma_s}\, \sigma_r\sigma_s, 
&\text{$s=r$},\\
\frac 1{\sigma_r\sigma_s} (0-\mu_r\mu_s-\mu_r\mu_s+\mu_r\mu_s), &
\text{$s>r+1$},\\[1em]
\frac 1{\sigma_r\sigma_s}(-p_s-\mu_r\mu_s-\mu_r\mu_s+\mu_r\mu_s), &
\text{$s=r+1$},\end{cases}\\
&=t\begin{cases} 1 & \text{$s=r$},\\
-\frac{\mu_r\mu_s}{\sigma_r\sigma_s} & \text{$s>r+1$},\\
-\frac{(p_s+\mu_r\mu_s)}{\sigma_r\sigma_s} & \text{$s=r+1$},
\end{cases}\end{align*}}
\noindent using the properties of the $Z^r_k$ noted at the beginning
of the proof.\CQFD 

\end{Proof}

We now study \eqref{item6zz} on a case-by-case basis.  
First, let $I^* = \emptyset$, 
that is, let $p_i = 1/m$, for $i = 1, \dots, m$. 
Then $\sigma_i^2 = 2p_i = 2/m$, for all $i \in \{1,2,\dots,m\}$.  
Hence, simply rescaling \eqref{item6zz} 
by $\sqrt{2/m}$ recovers the uniform result in \eqref{item6x}.

Next, consider the case where $p_{max}=p_j$, 
for {\it precisely one} $j\in \{1,\dots, m\}$.  
We then have 
$I^*=\{1,2,\dots, m\}\backslash \{j\}$.  
This forces us to set 
$0 = t_0 =  t_1 = \cdots = t_{j-1}$ 
and $t_j =  t_{j+1} = \cdots = t_{m-1} = t_m = 1$,  
in the maximizing term in \eqref{item6zz} .  
This leads to the following result.

\begin{corollary}\label{cor1}

If $p_{max}=p_j$ for precisely one $j\in \{1,\dots, m\}$, then

\begin{equation}\label{item16}
\frac{LI_n-p_{max}n}{\sqrt n} \Rightarrow 
-\frac1m\sum^{m-1}_{i=1} i\sigma_i\tilde{B}^i(1)
+\sum^{m-1}_{i=j} \sigma_i\tilde{B}^i(1),
\end{equation}

\noindent where the last term in \eqref{item16} is not present if $j=m$.

\end{corollary}

\begin{Rem} (i) Above, $(LI_n-p_{max}n)/\sqrt n$ 
converges to a centered normal random variable.
Intuitively, this result is not surprising 
since the longest increasing subsequence is, 
asymptotically, a string consisting 
primarily of the most frequently occurring 
letter, a string whose length is approximated 
by a binomial random variable with parameters $n$ 
and $p_{max}$. We show below that the 
variance of the limiting normal distribution is, in fact,
equal to $p_{max}(1-p_{max})$. \\
\noindent (ii) One could compute the variance of the 
right hand side of \eqref{item16} directly to verify that
it is in fact $p_{max}(1-p_{max})$.  However, the nature of 
the covariance structure of the Brownian motion
makes the calculation somewhat cumbersome.  Instead, 
we revisit the appoximation to our Brownian motion 
in the first term on the right hand side of \eqref{item16}. 
In doing this, we not only recover the variance
of the limiting distribution, but also see 
that our approximating functional does indeed 
take the form of the sum of a binomial random variable 
and of a term which converges to zero in probability. 
\end{Rem}

\noindent\begin{Proof} We have from the very 
definition of the approximation to Brownian motion that

\begin{align}\label{item17}
- \frac1m\sum^{m-1}_{i=1} i\sigma_i\hat B^i_n(1) 
  &= - \frac1m\sum^{m-1}_{i=1} i\sigma_i\ \left[\frac{S^i_n-\mu_i n}
{\sigma_i\sqrt n}\right]\nonumber\\
  &=   \frac{1}{\sqrt n} \left[- \frac1m\sum^{m-1}_{i=1} i S^i_n 
                                 + \frac{n}{m}\sum^{m-1}_{i=1} i\mu_i \right].
\end{align}

\noindent Recalling that 
$- \frac1m \sum^{m-1}_{i=1} i S^i_n = a^m_n - \frac{n}{m}$, and that 
$\sum^{m-1}_{i=1} i\mu_i = 1 - mp_m$, 
\eqref{item17} becomes

\begin{equation}\label{item18}
\frac{1}{\sqrt n}  \left[ \left(a^m_n - \frac{n}{m}\right) + 
\frac{n}{m}(1 - mp_m) \right]
= \frac{1}{\sqrt n} (a^m_n - np_m).
\end{equation}

Turning to the second term on the 
right hand side of \eqref{item16} and noting 
that for $1 \le j < k \le {m-1}$,
$\sum^k_{i=j} \mu_i = p_j - p_{k+1}$ 
and that $\sum^k_{i=j} S^i_r = a^j_r - a^{k+1}_r$,
for $1 \le r \le n$, we then have

\begin{align}\label{item19}
\sum^{m-1}_{i=j} \sigma_i\hat B^i_n(1) 
  &= \frac{1}{\sqrt n} \left[ \sum^{m-1}_{i=j} S^i_n  - n 
\sum^{m-1}_{i=j} \mu_i \right]\nonumber\\
  &= \frac{1}{\sqrt n} \left[ (a^j_n - a^m_n) - n(p_j-p_m) 
\right]\nonumber\\
  &= \frac{1}{\sqrt n} \left[ (a^j_n - np_j) - (a^m_n - np_m) \right].
\end{align}

We saw in \eqref{item14} that we 
could write $(LI_n - p_{max}n)/\sqrt{n}$,
as the sum of a functional 
approximating the Brownian motion and 
of an error term $E_n$ converging, to zero, in probability. 
In the present case, this expression simplifies to 

\begin{equation}\label{item20}
-\frac1m\sum^{m-1}_{i=1} i\sigma_i\hat B^i(1) +\sum^{m-1}_{i=j} 
\sigma_i\hat B^i(1) + E_n
= \frac{a^j_n - np_j}{\sqrt n} + E_n,
\end{equation}

\noindent using
\eqref{item17}--\eqref{item19}. 

Now $a^j_n$ is a binomial random variable with parameters 
$n$ and $p = p_j = p_{max}$.  By the Central
Limit Theorem and the converging together lemma, 
the right hand side of \eqref{item20} converges to a 
$N(0,p_{max}(1-p_{max}))$ distribution, while
by Donsker's Theorem, the left hand side converges to the 
Brownian functional obtained in \eqref{item16}.  Hence,  
$(LI_n - p_{max}n)/\sqrt n \Rightarrow N(0,p_{max}(1-p_{max}))$, 
as claimed.\CQFD 
\end{Proof}

Let us now study what happens when 
$p_{max} = p_j = p_k$, $1 \le j < k \le m$,
and $p_i < p_{max}$ otherwise, that is, 
when {\it precisely two letters} 
have the maximal probability.  
We then have 
$I^* = \{1,\dots,m\}\backslash\{j,k\}$.  
This requires that

\begin{align}
&0 = t_0 = t_1 = \cdots = t_{j-1},\nonumber\\
&t_j = t_{j+1} = \cdots = t_{k-1},\nonumber\\
&t_k = t_{k+1} = \cdots = t_m = 1.\nonumber
\end{align}

Hence,
{\allowdisplaybreaks
\begin{align}
\mathop{\max_{\stackrel{\scriptstyle 0 = t_0 \le t_1\le \cdots}{\le t_{m-1}
\le t_m = 1}}}\sum^{m-1}_{i=1} \sigma_i \tilde{B}^i(t_i)
&=\max_{0 \le t \le 1} \left[  \sum^{k-1}_{i=j} \sigma_i \tilde{B}^i(t) + 
\sum^{m-1}_{i=k} \sigma_i \tilde{B}^i(1) \right]\nonumber\\
&=\sum^{m-1}_{i=k} \sigma_i \tilde{B}^i(1) + \max_{0 \le t \le 1}   
\sum^{k-1}_{i=j} \sigma_i \tilde{B}^i(t)\nonumber.
\end{align}
}

Thus the limiting law is 
\begin{equation}\label{item20b}
- \frac1m\sum^{m-1}_{i=1} i\sigma_i \tilde{B}^i_n(1) 
+ \sum^{m-1}_{i=k} \sigma_i \tilde{B}^i(1) + \max_{0 \le t \le 1}
\sum^{k-1}_{i=j} \sigma_i \tilde{B}^i(t).
\end{equation}

To consolidate our analysis, we treat the general case for 
which $p_{max}$ occurs exactly $k$ times among 
$\{p_1,p_2,\dots,p_m\}$, where $2 \le k \le m-1$.  
%Note that we must have $1/m \le p_{max} \le 1/k$.  
Not only will we recover the natural analogues of 
\eqref{item20b}, but we will also express our 
results in terms of another functional of
Brownian motion which is more symmetric.  
Combining the $2 \le k \le m-1$ 
case at hand with the $k=1$ 
case previously examined, we have the following:\\

\begin{corollary}\label{cor2}
Let $p_{max} = p_{j_1} =p_{j_2} = \cdots \ = p_{j_k}$ for 
$1 \le j_1 < j_2 < \cdots < j_k \le m,$
for some $1 \le k \le m-1$, 
and let $p_i < p_{max}$, otherwise. Then

\begin{equation}\label{item25}
\frac{LI_n-p_{max}n}{\sqrt n} \Rightarrow 
\sqrt{p_{max} (1-p_{max}) }  
   \mathop{\max_{\stackrel{\scriptstyle 0 = t_0 \le t_1\le \cdots}
{\le t_{k-1} \le t_k = 1}}}
     \sum^{k}_{\ell=1} \left[ \tilde{B}^{\ell}(t_{\ell}) - 
\tilde{B}^{\ell}(t_{\ell-1}) \right],
\end{equation}

\noindent where the $k$-dimensional Brownian motion 
$(\tilde{B}^1(t),\tilde{B}^1(t),\dots,\tilde{B}^k(t))$
has the covariance matrix

\begin{equation}\label{item25a}t\begin{pmatrix}
1 & \rho & \rho & \cdots & \rho\\
\rho & 1 & \rho && \vdots\\
\vdots &\ddots &\ddots & \ddots & \vdots\\
\vdots &&\rho & 1 & \rho\\
\rho &\cdots &\cdots &\rho &1\end{pmatrix},
\end{equation}

\noindent with $\rho = - p_{max}/(1-p_{max})$. 

\end{corollary}

\noindent \begin{Proof} Let 
$p_{max} = p_{j_1} =p_{j_1} = \cdots \ = p_{j_k},$ 
with $1 \le j_1 < j_2 < \cdots < j_k \le m$
and  $2 \le k \le m-1$, {\it i.e.}, let 
$I^* = \{1,2,\dots,m\} \backslash \{j_1,j_2,\dots,j_k\}$.  
Set $j_0 = 1$ and $j_{k+1} = m$.  
Then \eqref{item14} becomes

{\allowdisplaybreaks
\begin{align}\label{item21}
\frac{LI_n-p_{max}n}{\sqrt n} &= 
- \frac1m  \sum^{m-1}_{i=1} i\sigma_i\hat B^i_n(1)\nonumber\\
&\qquad + \mathop{\max_{\stackrel{\scriptstyle 0 = t_0 \le t_1\le \cdots}
{\le t_{m-1}\le t_m = 1}}}_{t_i=t_{i-1},\ i\in I^*}\sum^{m-1}_{i=1} 
\sigma_i\hat B^i_n(t_i) + E_n\nonumber\\
&= - \frac1m  \sum^{m-1}_{i=1} i\sigma_i\hat B^i_n(1)\nonumber\\
&\qquad +\mathop{\max_{\stackrel{\scriptstyle 0 = t_{j_0} \le t_{j_1}\le \cdots}
{\le t_{j_k}\le t_{j_{k+1}} = 1}}}\sum^{k}_{\ell=0} 
\sum^{j_{\ell+1}-1}_{i={j_\ell}} \sigma_i\hat B^i_n(t_{j_l}) + E_n\nonumber\\
&= - \frac1m  \sum^{m-1}_{i=1} i\sigma_i\hat B^i_n(1)+\nonumber\\
&\qquad + \mathop{\max_{\stackrel{\scriptstyle 0 = t_{j_0} \le t_{j_1}\le \cdots}
{\le t_{j_k}\le t_{j_{k+1}} = 1}}}  \left[ \sum^{k-1}_{\ell=1} 
\sum^{j_{\ell+1}-1}_{i={j_\ell}} \sigma_i\hat B^i_n(t_{j_\ell})+ \sum^{m-1}_{i={j_k}} 
\sigma_i\hat B^i_n(1) \right]\nonumber + E_n\\
&= \left[ - \frac1m  \sum^{m-1}_{i=1} i\sigma_i\hat B^i_n(1)+ 
\sum^{m-1}_{i={j_k}} \sigma_i\hat B^i_n(1)\right]\nonumber\\
&\qquad  + \mathop{\max_{\stackrel{\scriptstyle 0 = t_{j_0} \le t_{j_1}\le \cdots}
{\le t_{j_k}\le t_{j_{k+1}} = 1}}}
 \sum^{k-1}_{\ell=1} \sum^{j_{\ell+1}-1}_{i={j_\ell}} \sigma_i\hat B^i_n(t_{j_\ell}) 
+ E_n.
\end{align}}

We immediately recognize the first term 
on the right hand side of \eqref{item21} 
as what we encountered for $k=1$.  Using 
the definition of the $\hat B^i_n$, 
\eqref{item21} can then be rewritten as

{\allowdisplaybreaks
\begin{align}\label{item22}
&
\frac{a^{j_k}_n - np_{max}}{\sqrt n}
             + \mathop{\max_{\stackrel{\scriptstyle 0 
= t_{j_0} \le t_{j_1}\le \cdots}
{\le t_{j_k}\le t_{j_{k+1}} = 1}}}
                \sum^{k-1}_{\ell=1} \sum^{j_{\ell+1}-1}_{i={j_\ell}} 
\sigma_i\hat B^i_n(t_{j_\ell}) + E_n\nonumber\\
&\qquad =
\frac{a^{j_k}_n - np_{max}}{\sqrt n}
               + \mathop{\max_{\stackrel{\scriptstyle 0 = t_{j_0} 
\le t_{j_1}\le \cdots}{\le t_{j_k}\le t_{j_{k+1}} = 1}}}                
\sum^{k-1}_{\ell=1} \sum^{j_{\ell+1}-1}_{i={j_\ell}} \sigma_i 
		    \left( \frac{S^i_{\left[nt_{j_\ell}\right]}-\mu_i n}
{\sigma_i\sqrt n} \right) + E_n\nonumber\\
&\qquad=
\frac{a^{j_k}_n - np_{max}}{\sqrt n}
             + \frac{1}{\sqrt n} \mathop{\max_{\stackrel{\scriptstyle 0 = t_{j_0} 
\le t_{j_1}\le \cdots}{\le t_{j_k}\le t_{j_{k+1}} = 1}}}
                \sum^{k-1}_{\ell=1}  
		    \left[ (a^{j_\ell}_{\left[nt_{j_\ell}\right]} - 
a^{j_{\ell+1}}_{\left[nt_{j_\ell}\right]}) - n(p_{j_\ell}-p_{j_{\ell+1}}) \right]\nonumber\\
		    &\qquad\qquad + E_n \nonumber\\
&\qquad=
\frac{a^{j_k}_n - np_{max}}{\sqrt n}
             + \frac{1}{\sqrt n}\! \mathop{\max_{\stackrel{\scriptstyle 0 = 
t_{j_0} \le t_{j_1}\le \cdots}{\le t_{j_k}\le t_{j_{k+1}} = 1}}}
                \sum^{k-1}_{\ell=1}\!   
		    \left[\!( a^{j_\ell}_{\left[nt_{j_\ell}\right]} - np_{max}\!) - 
(\!a^{j_{\ell+1}}_{\left[nt_{j_\ell}\right]} - np_{max}\!)\!\right] \nonumber\\
		    &\qquad\qquad\qquad\qquad\qquad + E_n.      
\end{align}}

Setting $a^{j_{k+1}}_n = n - \sum^{k}_{\ell=1} a^{j_\ell}_n$, 
we note that the random vector  
$(a^{j_1}_n,a^{j_2}_n,\dots,\\
a^{j_{k+1}}_n)$ 
follows a multinomial distribution with parameters 
$n$ and $(p_{max},p_{max},\\
\dots,p_{max},1-kp_{max})$.  
It is thus natural to introduce a new 
Brownian motion approximation as follows:

\begin{equation}\label{item23}
\check{B}^{\ell}_n(t) = 
      \frac{a^{j_\ell}_{\left[nt_{j_\ell}\right]} - np_{max}}{\sqrt{np_{max} 
(1-p_{max}) }},
      \quad \text{$1 \le \ell \le k$}.
\end{equation}

Substituting \eqref{item23} into \eqref{item22} gives

\begin{align}\label{item24}
&
\sqrt{p_{max} (1-p_{max}) } \left\{ \check{B}^k_n(1) +
   \mathop{\max_{\stackrel{\scriptstyle 0 = t_0 \le t_1\le \cdots}
{\le t_{k-1} \le t_k = 1}}}
     \sum^{k-1}_{\ell=1} \left[ \check{B}^{\ell}_n(t_{\ell}) - 
\check{B}^{\ell+1}_n(t_{\ell}) \right] \right\}  + E_n   \nonumber\\
&=
\sqrt{p_{max} (1-p_{max}) }  
   \mathop{\max_{\stackrel{\scriptstyle 0 = t_0 \le t_1\le \cdots}
{\le t_{k-1} \le t_k = 1}}}
     \sum^{k}_{\ell=1} \left[ \check{B}^{\ell}_n(t_{\ell}) - 
\check{B}^{\ell}_n(t_{\ell-1}) \right] + E_n.
\end{align}

By Donsker's Theorem, our approximations 
$(\check{B}^{1}_n(t),\check{B}^{2}_n(t),\dots,\check{B}^{k}_n(t))$ 
converges jointly to a $k$-dimensional Brownian motion 
$(\tilde{B}^1(t),\tilde{B}^2(t),\dots,\tilde{B}^k(t)).$  
This Brownian motion has the covariance structure

$$ t\begin{pmatrix}
1 & \rho & \rho & \cdots & \rho\\
\rho & 1 & \rho && \vdots\\
\vdots &\ddots &\ddots & \ddots & \vdots\\
\vdots &&\rho & 1 & \rho\\
\rho &\cdots &\cdots &\rho &1\end{pmatrix},$$

\noindent  where $\rho = -p_{max}/(1-p_{max})$, 
a fact which follows immediately from 
the covariance of the multinomial distribution, 
where the covariance of any two  
distinct $a^{j_\ell}_r$ is simply $-rp^2_{max}$, 
for $1 \le r \le n$.  This, together with 
our analysis of the unique $p_{max}$ case,
proves the corollary.  \CQFD 

\end{Proof}

\begin{Rem}\label{rem1}
The above results provide a 
Brownian functional equivalent to the
GUE result of Its, Tracy, and Widom 
\cite{ITW1} (described in detail in
the comments preceding Theorem \ref{thm2} and with a law 
given in Remark \ref{sumrmt}).  
Note that the limiting 
distribution in \eqref{item25} depends only 
on $k$ and $p_{max}$; neither the 
specific values of $j_1, j_2, \dots, j_k$ 
nor the remaining values of $p_i$ are 
material, a fact already noted in \cite{ITW1}.  
Also, it follows from generic results on Brownian functionals
that this limiting law has a density,
which in the uniform case is supported on the
positive real line, while supported on all of $\bbr$ 
in the non-uniform case.
\end{Rem}

We have already seen in \eqref{item6z} that 
the limiting distribution for the uniform case
has a nice representation as a functional 
of standard Brownian motion.  We now also 
express the limiting distribution in \eqref{item25} 
as a functional of standard Brownian motion.
Moreover, this new functional extends to
the uniform case, although its form 
is different from that of \eqref{item6z}.
This limiting random variable can 
be viewed as the sum of a normal one and 
of a maximal eigenvalue type one.

\begin{corollary}\label{cor3}
Let $p_{max} = p_{j_1} =p_{j_2} = \cdots \ = p_{j_k},$ 
for $1 \le j_1 < j_2 < \cdots < j_k \le m$,
and some $1 \le k \le m$, and let 
$p_i < p_{max}$, otherwise. Then

\begin{align}\label{item24b}
&\frac{LI_n-p_{max}n}{\sqrt n} \Rightarrow \sqrt{p_{max}} 
   \Bigl\{  \frac{\sqrt{1-kp_{max}}-1}{k}
     \sum^{k}_{j=1} B^{j}(1)   \nonumber\\
   &\qquad \qquad \qquad \qquad+   
    \mathop{\max_{\stackrel{\scriptstyle 0 = t_0 \le t_1\le \cdots}
{\le t_{k-1} \le t_k = 1}}}
    \sum^{k}_{\ell=1}  \left[ B^{\ell}(t_{\ell}) 
          - B^{\ell}(t_{\ell-1}) \right]   \Bigr\}.
\end{align}

\noindent where $(B^1(t),B^2(t),\dots,B^k(t))$ 
is a standard $k$-dimensional Brownian motion.
 
\end{corollary}

\noindent \begin{Proof} 
Let us first examine the non-uniform case
$1 \le k \le m-1$.  Recall that 
$\rho = -p_{max}/(1-p_{max})$.
Now the covariance matrix in \eqref{item25a}
has eigenvalues $\lambda_1=1-\rho = 1/(1-p_{max})$ 
of multiplicity $k-1$ and 
$\lambda_2=1+(k-1)\rho = (1-kp_{max})/(1-p_{max}) < \lambda_1$ 
of multiplicity $1$.  
From the symmetries of the covariance matrix,
it is not hard to see that we can write 
each Brownian motion $\tilde{B}^i(t)$ as a 
linear combination of standard Brownian 
motions $(B^1(t),\dots,B^k(t))$
as follows:

\begin{align}\label{item26}
\tilde{B}^i(t) = \beta B^i(t) + \eta \sum^k_{j=1,j \ne i}B^j(t), && \text{$i=1,\dots,k$},
\end{align}

\noindent where 
\begin{align}\label{item27}
\beta = \frac{(k-1)\sqrt{\lambda_1} + \sqrt{\lambda_2}}{k},  && 
\eta = \frac{-\sqrt{\lambda_1} + \sqrt{\lambda_2}}{k}.
\end{align}

Substituting \eqref{item26} and 
\eqref{item27} into \eqref{item25},
and noting that 
$\beta - \eta = \sqrt{\lambda_1} = 1/\sqrt{1-p_{max}}$,
we find that

{\allowdisplaybreaks
\begin{align}\label{item28}
&\sqrt{p_{max} (1-p_{max}) }  
   \mathop{\max_{\stackrel{\scriptstyle 0 = t_0 \le t_1\le \cdots}{\le t_{k-1} \le t_k = 1}}}
     \sum^{k}_{\ell=1} \left[ \tilde{B}^{\ell}(t_{\ell}) - \tilde{B}^{\ell}(t_{\ell-1}) \right] \nonumber\\
&=\sqrt{p_{max} (1-p_{max}) } 
   \mathop{\max_{\stackrel{\scriptstyle 0 = t_0 \le t_1\le \cdots}{\le t_{k-1} \le t_k = 1}}}
        \sum^{k}_{\ell=1} \Bigl\{ \beta \left[ B^{\ell}(t_{\ell}) - B^{\ell}(t_{\ell-1}) \right]\nonumber\\
                                 &\qquad\qquad+\eta  \sum^{k}_{j=1,j\ne \ell} 
\left[ B^{j}(t_{\ell}) - B^{j}(t_{\ell-1}) \right]  \Bigr\}\nonumber\\
&=\sqrt{p_{max} (1-p_{max}) }  
   \mathop{\max_{\stackrel{\scriptstyle 0 = t_0 \le t_1\le \cdots}{\le t_{k-1} \le t_k = 1}}}
        \sum^{k}_{\ell=1} \Bigl\{ (\beta-\eta) \left[ B^{\ell}(t_{\ell}) - B^{\ell}(t_{\ell-1}) \right] \nonumber\\
                                      &\qquad\qquad+\eta \sum^{k}_{j=1} \left[ B^{j}(t_{\ell}) - B^{j}(t_{\ell-1}) \right]  \Bigr\}\nonumber\\
&=\sqrt{p_{max} (1-p_{max}) }  
   \mathop{\max_{\stackrel{\scriptstyle 0 = t_0 \le t_1\le \cdots}{\le t_{k-1} \le t_k = 1}}}
        \Bigl\{  \sum^{k}_{\ell=1}  (\beta-\eta) \left[ B^{\ell}(t_{\ell}) - B^{\ell}(t_{\ell-1}) \right] \nonumber\\
               &\qquad\qquad+ \eta   \sum^{k}_{\ell=1}  
\sum^{k}_{j=1} \left[ B^{j}(t_{\ell}) - B^{j}(t_{\ell-1}) \right]  \Bigr\}\nonumber\\
&=\sqrt{p_{max} (1-p_{max}) } \Bigl\{  \eta  \sum^{k}_{j=1} B^{j}(1)   \nonumber\\
  &\qquad\qquad+ (\beta-\eta)  \mathop{\max_{\stackrel{\scriptstyle 0 = t_0 \le t_1\le \cdots}{\le t_{k-1} \le t_k = 1}}}
                    \sum^{k}_{\ell=1}  \left[ B^{\ell}(t_{\ell}) - B^{\ell}(t_{\ell-1}) \right]   \Bigr\}\nonumber\\
&= \sqrt{p_{max}} 
   \Bigl\{  \frac{\sqrt{1-kp_{max}}-1}{k}
     \sum^{k}_{j=1} B^{j}(1)   \nonumber\\
   &\qquad \qquad+   
    \mathop{\max_{\stackrel{\scriptstyle 0 = t_0 \le t_1\le \cdots}{\le t_{k-1} \le t_k = 1}}}
    \sum^{k}_{\ell=1}  \left[ B^{\ell}(t_{\ell}) 
          - B^{\ell}(t_{\ell-1}) \right]   \Bigr\}.
\end{align}}

To complete the proof, we now examine
the uniform case $k=m$, where necessarily
$p_{max} = 1/m$.
Now we saw in Proposition
\ref{prop0} that

\begin{equation}\label{item28a}
\frac{LI_n-n/m}{\sqrt{n}} \Rightarrow 
\sqrt{\frac{2}{m}}
\Bigl\{ -\frac1m\sum^{m-1}_{i=1} i \tilde{B}^i(1)+
\max_{\stackrel{\scriptstyle 0\le t_1\le\cdots}{\le t_{m-1}\le 1}}
\sum^{m-1}_{i=1} \tilde{B}^i(t_i)\Bigr\},
\end{equation}

\noindent where the $(m-1)$-dimensional Brownian motion 
$(\tilde{B}^1(t),\dots,\tilde{B}^{m-1}(t))$ 
had a tridiagonal covariance matrix given by 
\eqref{mattri}. Now we can derive this Brownian motion
from a standard $m$-dimensional
Brownian motion
$(B^1(t),\dots,B^{m}(t))$
via the a.s.~transformations

\begin{align*}
\tilde{B}^i(t)     &= \frac{1}{\sqrt{2}}(B^{i}(t)   - 
B^{i+1}(t)), \qquad \text{$1\le i\le m-1$}.
\end{align*}

It is easily verified that the 
Brownian motion $(\tilde{B}^1(t),\dots,\tilde{B}^{m-1}(t))$ 
so obtained does indeed have
the covariance structure given by \eqref{mattri}.
Substituting these independent Brownian
motions into \eqref{item28}, we obtain
the following a.s~equalities:

{\allowdisplaybreaks
\begin{align}\label{item28b}
\frac{LI_n-n/m}{\sqrt{n}} 
&\Rightarrow 
  \sqrt{\frac{2}{m}}
  \Bigl\{ -\frac1m\sum^{m-1}_{i=1} i \tilde{B}^i(1)+
  \max_{\stackrel{\scriptstyle 0\le t_1\le\cdots}{\le t_{m-1}\le 1}}
  \sum^{m-1}_{i=1} \tilde{B}^i(t_i)\Bigr\}\nonumber\\
&= \sqrt{\frac{1}{m}}
  \Bigl\{ -\frac1m\sum^{m-1}_{i=1} i [B^i(1)-B^{i+1}(1)] \nonumber\\
  &\qquad + \max_{\stackrel{\scriptstyle 0\le t_1\le\cdots}{\le t_{m-1}\le 1}}
  \sum^{m-1}_{i=1} [B^i(t_i)-B^{i+1}(t_i)]\Bigr\}\nonumber\\
&= \sqrt{\frac{1}{m}}
  \Bigl\{ -\frac1m\sum^{m}_{i=1} B^i(1) + B^m(1) \nonumber\\
  &\qquad +\max_{\stackrel{\scriptstyle 0\le t_1\le\cdots}{\le t_{m-1}\le 1}}
  \sum^{m}_{i=1} [B^i(t_i)-B^{i}(t_{i-1})] - B^m(1)  \Bigr\}\nonumber\\
&= \sqrt{\frac{1}{m}}
  \Bigl\{ -\frac1m\sum^{m}_{i=1} B^i(1) 
  +\max_{\stackrel{\scriptstyle 0\le t_1\le\cdots}{\le t_{m-1}\le 1}}
  \sum^{m}_{i=1} [B^i(t_i)-B^{i}(t_{i-1})]  \Bigr\},
\end{align}}

\noindent which we recognize as 
\eqref{item24b}, with $k=m$ and
$p_{max} = 1/m$.\CQFD 
\end{Proof}

We have already seen several 
representations for the limiting law 
in the uniform case. Yet one
more pleasing functional
for the limiting distribution of $LI_n$
is described in the following

\begin{theorem}\label{thm3}
Let $p_{max} = p_{1} = p_{2} = \cdots \ = p_{m} = 1/m$.
Then 

\begin{equation}\label{item29}
\frac{LI_n - n/m}{\sqrt{n}} \Rightarrow
\frac{1}{\sqrt{m}} \mathop{\max_{\stackrel{\scriptstyle 0 = 
t_0 \le t_1\le \cdots}{\le t_{m-1} \le t_m = 1}}}
\sum^{m}_{i=1}  \left[ \tilde{B}^{i}(t_{i}) - \tilde{B}^{i}(t_{i-1}) \right]
:= \frac{\tilde{H}_m}{\sqrt{m}},
\end{equation}

\noindent where $(\tilde{B}^1(t),\tilde{B}^2(t),\dots,\tilde{B}^{m}(t))$ is an
$m$-dimensional Brownian motion having covariance
matrix \eqref{item25a}, with $\rho = -1/(m-1)$. 
(This Brownian motion satisfies
$\sum_{i=1}^{m} \tilde{B}^{i}(t) = 0$, 
for all $0 \le t \le 1$.)

\end{theorem}

\noindent \begin{Proof}  We show that the 
functional being maximized in \eqref{item29} 
has the same covariance structure as the 
functional being maximized in \eqref{item6z}, 
a result which we restate as:

\begin{equation}\label{item30}
\frac{LI_n-n/m}{\sqrt{n}} \Rightarrow 
\frac{1}{\sqrt{m}}  \max_{\stackrel{\scriptstyle 0\le t_1\le\cdots}
{\le t_{m-1}\le t_m = 1}}
                   \sum^{m-1}_{i=1} \left[ \beta_i B^i(t_{i+1}) - 
\eta_i B^i(t_i) \right],
\end{equation}

\noindent where $\beta_i =\sqrt{i/(i+1)}$ 
and $\eta_i = \sqrt{(i+1)/i}$.
From this it will immediately 
follow that the maxima, over all 
$0 \le t_1 \le t_2 \le \dots \le t_{m-1} \le 1$,
in both expressions have the 
same law, clinching the proof.

Given that the zero-sum condition 
on the Brownian motion is in force 
in \eqref{item29}, it is natural 
to rewrite \eqref{item29} as

\begin{align}\label{item31}
&\frac{LI_n-n/m}{\sqrt n} \Rightarrow  
  \frac{1}{\sqrt{m}} \mathop{\max_{\stackrel{\scriptstyle 0 = t_0 \le t_1\le \cdots}
{\le t_{m-1} \le t_m = 1}}}
                    \sum^{m}_{i=1}  \left[ \tilde{B}^{i}(t_{i}) - \tilde{B}^{i}(t_{i-1}) \right],
\end{align}

\noindent where  
$(\tilde{B}^1(t),\tilde{B}^2(t),\dots,\tilde{B}^m(t))$ 
is an $m$-dimensional Brownian motion
with a permutation-invariant 
covariance matrix described by  

\begin{align*}
\mbox{Cov} (\tilde{B}^i(t),\tilde{B}^j(t)) &= \left(\frac{m-1}{m}\right)\frac{-t}{m-1}\\
 &= -\frac{t}{m}, \quad i \ne j,
\end{align*}

\noindent while

\begin{align*}
\mbox{Var}  \tilde{B}^i(t) &= \frac{m-1}{m}t.
\end{align*}

Let $t = (t_1, t_2, \dots , t_{m-1})$ 
be a fixed collection of $t_i$ from the Weyl chamber 
$T = \{(t_1, t_2, \dots , t_{m-1}):0 \le t_1 \le t_2 \le \dots \le t_{m-1} \le 1\}$.  
Setting 

\begin{equation}\label{item31a}
X_t = \sum^{m}_{i=1}  \left[ \tilde{B}^{i}(t_{i}) - \tilde{B}^{i}(t_{i-1}) \right],
\end{equation}

\noindent we then have 

\begin{align}\label{item32}
\mbox{Cov} (X_t,X_s) 
  &= \frac{m-1}{m} \sum_{1 \le i,j \le m}  \mbox{Cov} (\tilde{B}^{i}(t_{i}) - 
\tilde{B}^{i}(t_{i-1}), \tilde{B}^{i}(s_{i}) - \tilde{B}^{i}(s_{i-1})) \nonumber\\
  &= \frac{m-1}{m} \sum^{m}_{i=1} \left[t_{i} \wedge s_{i} - t_{i} 
\wedge s_{i-1} - t_{i-1} \wedge s_{i} + t_{i-1} \wedge s_{i-1}     \right] \nonumber\\
  & \quad - \frac{1}{m} \sum_{i \ne j} \left[t_{i} \wedge s_{j} - t_{i} 
\wedge s_{j-1} - t_{i-1} \wedge s_{j} + t_{i-1} \wedge s_{j-1}     \right].
\end{align}

\noindent  We can rewrite \eqref{item32} 
in an especially clear way by setting
 $T_1 = [0,t_1]$ and $T_i = (t_i,t_{i+1}]$, $i = 2, \dots, m$,
and similarly $S_1 = [0,s_1]$ and 
$S_i = (s_i,s_{i+1}]$, $i = 2, \dots, m$.  
Letting  $Leb$ denote the Lebesgue measure on $[0,1]$, 
a case-by-case analysis of the relative positions 
of $t_i, t_{i-1}, s_i$, and $s_{i-1}$ quickly yields that

\begin{align}\label{item33}
\mbox{Cov} (X_t,X_s)
 &= \frac{m-1}{m} \sum^{m}_{i=1} Leb(T_i \cap S_i) - 
\frac{1}{m} \sum_{i \ne j} Leb(T_i \cap S_j) \nonumber\\
 & = \frac{m-1}{m} \sum^{m}_{i=1} Leb(T_i \cap S_i) - 
\frac{1}{m} \left[ 1 -\sum^{m}_{i=1} Leb(T_i \cap S_i) \right] \nonumber\\
 &= - \frac{1}{m} + \sum^{m}_{i=1} Leb(T_i \cap S_i).
\end{align}

\noindent From \eqref{item33} we clearly have 
Var$X_t = (m-1)/m$, for all $t \in T$.  
To complete the proof, we now show that 

\begin{equation}\label{item34}
Y_t = \sum^{m-1}_{i=1} \left[  \beta_i B^i(t_{i+1}) -  \eta_i B^i(t_i) \right],
\end{equation}

\noindent has the same covariance structure as $X_t$, 
where $\beta_i =\sqrt{i/(i+1)}$ 
and $\eta_i = \sqrt{(i+1)/i}$.
Using the independence of the components 
of the Brownian motion, we also have

\begin{align}\label{item35}
\mbox{Cov} (Y_t,Y_s) 
  &= \sum^{m-1}_{i=1} \mbox{Cov} \left(\beta_i B^i(t_{i+1}) - \eta_i B^i(t_i), 
\beta_i B^i(s_{i+1}) -  \eta_i B^i(s_i)\right) \nonumber\\
  &= \sum^{m-1}_{i=1} \left[ \frac{i}{i+1} (t_{i+1} \wedge s_{i+1}) - t_{i+1} 
\wedge s_{i} - t_{i} \wedge s_{i+1} + \frac{i+1}{i} (t_{i} \wedge s_{i})  \right] 
\nonumber\\
  &= \sum^{m}_{i=1}  \frac{i-1}{i} t_{i} \wedge s_{i}  - \sum^{m-1}_{i=1}  
\left[ t_{i+1} \wedge s_{i} - t_{i} \wedge s_{i+1}  -  \frac{i+1}{i} t_{i} 
\wedge s_{i} \right] \nonumber\\
  &= \frac{m-1}{m} - \sum^{m-1}_{i=1}  \left[t_{i+1} \wedge s_{i} - t_{i} 
\wedge s_{i+1}  -  2  (t_{i} \wedge s_{i}) \right].
\end{align}

As before, a simple case-by-case analysis of the summands 
in \eqref{item35} reveals that 

\begin{align}\label{item36}
\mbox{Cov} (Y_t,Y_s) 
  &= \frac{m-1}{m} - \left[ 1 - \sum^{m}_{i=1} Leb(T_i \cap S_i) \right] \nonumber\\
  &= - \frac{1}{m} + \sum^{m}_{i=1} Leb(T_i \cap S_i),
\end{align}

\noindent completing the proof. \CQFD 

\end{Proof}

\section{Large-$m$ Asymptotics and Related Results}

With the covariance structure of $X_t$ now in hand, we can
compute the $L^2$-distance between any $X_t$ and $X_s$:

\begin{align}\label{item37}
\bbe(X_t - X_s)^2 &= \textvar X_t + \textvar X_s - 2\textcov (X_t,X_s) 
\nonumber\\
                    &= 2(1-1/m) - 2\left[-1/m + \sum^{m}_{i=1} 
Leb(T_i \cap S_i) \right]  \nonumber\\
		    &= 2\left[1- \sum^{m}_{i=1} Leb(T_i \cap S_i) \right].
\end{align}

\noindent Such a metric is useful, for instance, 
in applying Dudley's Entropy Bound to show
that

\begin{equation}\label{item37b}
\bbe \left(\mathop{\max_{\stackrel{\scriptstyle 0 
= t_0 \le t_1\le \cdots}{\le t_{m-1} \le t_m = 1}}} X_t\right) 
\le K\sqrt{m-1}, \nonumber
\end{equation}
for some constant $K$ not depending on $m$.

We can now more clearly see the similarities between the 
functional $D_m$ of Glynn and Whitt in \eqref{item6za} 
and that of \eqref{item6z}, which we have shown to have the 
same law as $\tilde{H}_m$ in \eqref{item29}.  Indeed, the only 
difference between the functionals is simply that in \eqref{item6za} 
the Brownian motions are independent, 
while in \eqref{item29} they are subject to the zero-sum constraint. 
Gravner, Tracy, and Widom \cite{GTW}
have already remarked that random words could
be studied via such Brownian functionals.
In fact, a restatement of
Corollary \ref{cor3} shows that, in law, 
$D_m$ and $\tilde{H}_m$ differ by a centered normal random variable, 
as indicated by the next theorem and corollary.
This, in turn, will allow us to clearly state asymptotic
results for $\tilde{H}_m$ from the known corresponding results
for $D_m$.

\begin{theorem}\label{thm4} 

Let $$H_m = \sqrt{2}\left\{-\frac1m\sum^{m-1}_{i=1} i \tilde{B}^i(1)+ 
\max_{\stackrel{\scriptstyle 0\le t_1\le\cdots}{\le t_{m-1}\le 1}} 
\sum^{m-1}_{i=1} \tilde{B}^i(t_i)\right\},$$

\noindent $m \ge 2$, and let $\tilde{H_1} \equiv 0$ a.s.,
where $(\tilde{B}^1(t),\dots,\tilde{B}^{m-1}(t))$ is an 
$(m-1)$-dimensional Brownian
motion with tridiagonal covariance matrix given by \eqref{mattri}.
Let $$D_m = \max_{\stackrel{\scriptstyle 0 = t_0 \le t_1\le\cdots}
{\le t_{m-1}\le t_m  = 1}}
\sum^m_{i=1} \left[ B^i(t_i) - B^i(t_{i-1}) \right],$$

\noindent where $(B^1(t),\dots,B^{m}(t))$ 
is a standard $m$-dimensional Brownian
motion.  Then $D_m = Z_m + H_m$ a.s., 
where $Z_m$ is a centered normal 
random variable with variance $1/m$, 
and in fact is given by
$Z_m = (1/m)\sum_{i=1}^m B^i(1)$.
\end{theorem}

\noindent \begin{Proof} 
The $m=1$ case is trivial.  For $m \ge 2$,
reformulating the proof of
Corollary \ref{cor3}, 
for the uniform case,
in terms of the
functionals $H_m$ and $D_m$
shows that

\begin{align*}
\frac{H_m}{\sqrt{m}} &= \frac{1}{\sqrt{m}}
\left(-\frac1m\sum_{i=1}^m B^i(1) + D_m\right)\\
            &= \frac{1}{\sqrt{m}}(-Z_m + D_m),
\end{align*}

\noindent almost surely, and hence
$D_m = Z_m + H_m$ a.s.\CQFD
\end{Proof}

Recalling the definition of $\tilde{H}_m$ from Theorem \ref{thm3}:
$$\tilde{H}_m := \mathop{\max_{\stackrel{\scriptstyle 0 = t_0 
\le t_1\le \cdots}{\le t_{m-1} \le t_m = 1}}}
\sum^{m}_{i=1}  \left[ \tilde{B}^{i}(t_{i}) - \tilde{B}^{i}(t_{i-1}) \right],$$
where $(\tilde{B}^{1}(t),\tilde{B}^{2}(t),\dots,\tilde{B}^{m}(t))$
is an $m$-dimensional Brownian motion having covariance matrix
\eqref{item25a}, with $\rho = -1/(m-1)$, 
{\it i.e.}, $\sum_{i=1}^{m} \tilde{B}^{i}(t) = 0$, 
for all $0 \le t \le 1$, we then have

\begin{corollary} \label{cor4}
For each $m \ge 1$, 
$\tilde{H}_m \stackrel{d}{=} D_m + Z_m$,
where $d$ denotes equality in distribution.
\end{corollary}

\noindent \begin{Proof} Proposition \ref{prop0} asserts that

$$\frac{LI_n-n/m}{\sqrt{n}} 
\Rightarrow \frac{H_m}{\sqrt{m}},$$ 

\noindent as $n \rightarrow \infty$, while by Theorem \ref{thm3}
$$\frac{LI_n-n/m}{\sqrt{n}} \Rightarrow \frac{\tilde{H}_m}{\sqrt{m}},$$ 

\noindent as $n \rightarrow \infty$ as well.  The conclusion
follows from the previous theorem.\CQFD 
\end{Proof}

This relationship between 
$\tilde{H}_m$ (resp.,$H_m$) and $D_m$ 
allows us to further express the limiting 
distribution in a rather compact form.

\begin{proposition}\label{prop1}
Let $p_{max} = p_{j_1} =p_{j_2} = \cdots \ = p_{j_k},$ 
for $1 \le j_1 < j_2 < \cdots < j_k \le m,$
and some $1 \le k \le m$, and let $p_i < p_{max}$, otherwise. Then

\begin{align*}
\frac{LI_n-p_{max}n}{\sqrt{n}} 
  &\Rightarrow \sqrt{p_{max}}\{\sqrt{1 - kp_{max}}Z_k + H_k\}\\
  &\stackrel{d}{=} \sqrt{p_{max}}\{\sqrt{1 - kp_{max}}Z_k + \tilde{H}_k\}.
\end{align*}

\end{proposition}

\noindent \begin{Proof} For $k=m$, we have 
$p_{max} = 1/m$, and thus simply
recover the limiting distribution 
$H_m/\sqrt{m} \stackrel{d}{=} \tilde{H}_m/\sqrt{m}$
of the uniform case.

For $1 \le k \le m-1$, we saw in Corollary \ref{cor3} 
that we could write the limiting law of 
$(LI_n-p_{max}n)/{\sqrt{n}}$ as

\begin{align}\label{item37db}
&\sqrt{p_{max}} 
   \Bigl\{  \frac{\sqrt{1-kp_{max}}-1}{k}
     \sum^{k}_{j=1} B^{j}(1)   \nonumber\\
   &\qquad \qquad+   
    \mathop{\max_{\stackrel{\scriptstyle 0 = t_0 \le t_1\le \cdots}
{\le t_{k-1} \le t_k = 1}}}
    \sum^{k}_{\ell=1}  \left[ B^{\ell}(t_{\ell}) 
          - B^{\ell}(t_{\ell-1}) \right]   \Bigr\},
\end{align}

\noindent where $(B^1(t),B^2(t),\dots,B^k(t))$ 
is a standard $k$-dimensional Brownian motion.
But, recalling the definitions of $D_k$ and
$Z_k$, and the fact that 
$D_k = Z_k + H_k$ a.s.,
\eqref{item37db} becomes

{\allowdisplaybreaks
\begin{align}\label{item37dc}
& \sqrt{p_{max}} 
   \left\{  \frac{\sqrt{1-kp_{\max}}-1}{k}(kZ_k) + D_k\right\}\nonumber\\
&= \sqrt{p_{max}} 
   \left\{  \left(\sqrt{1-kp_{max}}-1\right)Z_k + (Z_k + H_k)\right\}
\nonumber\\
&= \sqrt{p_{max}} 
   \left\{   \sqrt{1-kp_{max}} Z_k + H_k\right\}\nonumber\\
&\stackrel{d}{=} \sqrt{p_{max}} 
   \left\{   \sqrt{1-kp_{max}} Z_k + \tilde{H}_k\right\} .
\end{align}}\CQFD
\end{Proof}

\begin{Rem}
One might also write the limiting law of
Proposition \ref{prop1} in terms of the 
functional $D_k$.  Indeed, we have
$$\frac{LI_n-p_{max}n}{\sqrt{p_{max}n}}   
\Rightarrow (\sqrt{1 - kp_{max}}-1)Z_k + D_k,$$
so that the limiting law is expressed as
the sum of a centered normal random variable
and of the maximal eigenvalue of a $k\times k$ element of the GUE.
\end{Rem}

The behavior of $D_m$ has been well-studied.  
In particular, it is known 
that $D_m/\sqrt{m} \rightarrow 2$ a.s.~and 
in $L^1$, as $m \rightarrow \infty$
(see \cite{Ba,GW,HMO,OY1,OY2,S}),
and that $(D_m-2\sqrt{m})m^{1/6} \Rightarrow F_2$,
as $m \rightarrow \infty$, 
where $F_2$ is the
Tracy-Widom distribution
(see \cite{Ba,GTW,TW,TW2}).  
From these results, the
asymptotics of $H_m$ follows.

\begin{theorem}\label{thm5} We have that 
$$\frac{H_m}{\sqrt{m}} \rightarrow 2$$ 
\noindent a.s.~and in $L^1$, as $m \rightarrow \infty$. 
Moreover,

\begin{equation}\label{item37dd}
\left(\frac{H_m}{\sqrt{m}} - 2\right)m^{2/3} \Rightarrow F_2,
\end{equation}

\noindent where $F_2$ is the Tracy-Widom distribution.
The same statements hold for $\tilde{H}_m$ 
in place of $H_m$.
\end{theorem}

\noindent \begin{Proof} From Theorem \ref{thm4} 
we have $D_m = Z_m + H_m$ a.s.,
where $Z_m = (1/m)\sum_{i=1}^m B^i(1)$.   

Clearly, $Z_m \rightarrow 0$ a.s.~and in $L^1$.
Thus, a.s.~and in $L^1$,
$$\lim_{m \rightarrow \infty} \frac{H_m}{\sqrt{m}} 
= \lim_{m \rightarrow \infty} \frac{D_m}{\sqrt{m}} .$$

Since this last limit is $2$, and since, for 
each $m \ge 1$, $H_m \stackrel{d}{=} \tilde{H}_m$,
it also follows that

$$\lim_{m \rightarrow \infty} 
\bbe \left| \frac{\tilde{H}_m}{\sqrt{m}}  - 2 \right| = 0.$$

We are thus left with proving the a.s.~covergence
to $2$ of $\tilde{H}_m/\sqrt{m}$. Since the variance 
of the functional being maximized in
the definition of $\tilde{H}_m$ equals 
$1-1/m$, the Gaussian concentration 
inequality then implies that 
$$\bbp(|\tilde{H}_m - \bbe \tilde{H}_m| > h) \le 2e^{\frac{-h^2}{2(1-\frac1m)}} 
< 2e^{\frac{-h^2}{2}}$$
for all $h>0$.  Then since $\bbe \tilde{H}_m/\sqrt{m} \rightarrow 2$ as 
$m \rightarrow \infty$ 
we have for $m$ large enough that 

\begin{align*}
\bbp \left(\left|\frac{	\tilde{H}_m}{\sqrt{m}} - 2\right| > h\right) 
&\le  \bbp \left(\left| \tilde{H}_m - \bbe \tilde{H}_m\right|  > 
\sqrt{m}\left(h - \left|\frac{\bbe \tilde{H}_m}{\sqrt{m}} - 2\right|\right)\right)\\
&\le \bbp \left(\left|\tilde{H}_m - \bbe \tilde{H}_m\right|  > \frac{\sqrt{m}h}{2}\right)\\
&< 2e^{\frac{-mh^2}{8}}
\end{align*}

\noindent This concentration result implies that 
$$\sum_{m=1}^{\infty}\bbp \left(\left|\frac{\tilde{H}_m}{\sqrt{m}} - 2\right| > h\right)  
\le \sum_{m=1}^{\infty}2e^{\frac{-mh^2}{8}} < \infty,$$
and the Borel-Cantelli lemma allows us to conclude
the proof of a.s.~convergence.

Turning to the limiting law, we know (\cite{Ba,GTW}) 
that $D_m$ has the
same distribution as the largest eigenvalue of the
$m \times m$ GUE.  Then the fundamental 
random matrix theory result of
Tracy and Widom \cite{TW2}
implies that

$$\left(\frac{D_m}{\sqrt{m}} - 2\right)m^{2/3} \Rightarrow F_2.$$

Since, moreover, $D_m = Z_m + H_m$, 
and since $Z_m$ has variance $1/m$, $Z_mm^{1/6} \Rightarrow 0$,
and so

\begin{align*}
\left(\frac{H_m}{\sqrt{m}} - 2\right)m^{2/3} 
                &= \left(\frac{D_m}{\sqrt{m}} - 2\right)m^{2/3} - 
Z_m m^{1/6} \Rightarrow F_2.
\end{align*}
Finally, $H_m \stackrel{d}{=} \tilde{H}_m$, 
and the same result
holds for $\tilde{H}_m$ in place of 
$H_m$.\CQFD

\end{Proof}

\begin{Rem}
\noindent (i) In the conclusion to \cite{TW}, Tracy and Widom
already derived \eqref{item37dd} by applying a scaling argument
to the limiting distribution of the uniform
alphabet case.  In our case we can moreover assert 
that a.s.~and in the mean, 
$$\lim_{k\to +\infty}\lim_{n\to +\infty}\frac{LI_n-p_{max}n}{\sqrt{kp_{max}n}} = 2,$$  
and that 

$$\left(\frac{LI_n-p_{max}n}{\sqrt{kp_{max}n}} - 2\right)k^{2/3} \Rightarrow F_2,$$
where the weak limit is first taken over $n$ and then over $k$.    

\noindent (ii)
Using scaling, subadditivity, and concentration
arguments found in
Hambly, Martin, and O'Connell \cite{HMO} and in O'Connell
and Yor \cite{OY1}, one could prove directly that
$\tilde{H}_m/\sqrt{m} \rightarrow 2$ a.s.
This could be accomplished by studying,
as do these authors, a process
version of $\tilde{H}_m$, {\it i.e.},

\begin{align*}
\tilde{H}_m(\varepsilon) :=
\mathop{\max_{\stackrel{\scriptstyle 0 = t_0 \le t_1\le \cdots}
{\le t_{m-1} \le t_m =
 \varepsilon}}}
\sum^{m}_{i=1} \left[ \tilde{B}^{i}(t_{i}) - \tilde{B}^{i}(t_{i-1}) \right],
\end{align*}

\noindent for $\varepsilon > 0$.  
With obvious notations, for all 
$\varepsilon > 0$ and $m \ge 1$,
$D_m(\varepsilon) = Z(\varepsilon) + H_m(\varepsilon)$, a.s.,
where $Z(\varepsilon) = (1/m)\sum_{i=1}^m B^i(\varepsilon)$.

\end{Rem}

To see in further detail how $D_m$ and $\tilde{H}_m$ are related, 
first note that $D_m \le D_{m+1}$ a.s. 
for $m \ge 1$, since $D_m$ can simply 
be obtained by restricting the right-most 
parameter $t_m$ to be $1$ in the definition of $D_{m+1}$. 
We now show a stochastic domination result between
$D_m$ and $\tilde{H}_m$. 

Recall that a random variable $X$ 
is said to {\it stochastically dominate} 
another random variable $Y$
$({\it i.e.}, X \ge_{st} Y)$ if for all 
$x \in \bbr$ we have $\bbp(X \ge x) \ge \bbp(Y \ge x)$.

\begin{proposition}\label{prop3}
$\tilde{H}_m \ge_{st} \sqrt{(1-1/m)}D_m$, for $m \ge 1$. The same
statement holds for $H_m$ in place of $\tilde{H}_m$.
\end{proposition}

\noindent \begin{Proof}  Since the $m=1$ case is trivial, let $m \ge 2$.
We saw in \eqref{item33} that the functional $X_t$ being maximized 
in the definition of $\tilde{H}_m$ had a 
covariance structure given by 
$\mbox{Cov} (X_t,X_s) = - 1/m + \sum^{m}_{i=1} Leb(T_i \cap S_i)$.  
A similar argument shows that the functional 
$U_t = \sum^m_{i=1} \left[ B^i(t_i) - B^i(t_{i-1}) \right]$ 
which is being maximized in the definition of $D_m$ has a 
covariance structure given by 
$\mbox{Cov} (U_t,U_s) = \sum^{m}_{i=1} Leb(T_i \cap S_i)$.  
Therefore, 
$$\mbox{Var}(\sqrt{(1-1/m)} U_t) = \mbox{Var}X_t = 1-1/m,$$ 
and

\begin{align*}
\mbox{Cov} (\sqrt{(1-1/m)} U_t,\sqrt{(1-1/m)} U_s) &= 
(1-1/m) \sum^{m}_{i=1} Leb(T_i \cap S_i)\\ 
&\ge \mbox{Cov} (X_t,X_s).
\end{align*}

By Slepian's Lemma we conclude that $\tilde{H}_m \ge_{st} \sqrt{(1-1/m)}D_m$.  
The final assertion follows from the equality in 
law between $\tilde{H}_m$ and $H_m$.\CQFD

\end{Proof} 

\begin{Rem}  Note that 

\begin{align*}
\bbe (X_t - X_s)^2 &= \bbe (\sqrt{(1-1/m)} U_t - \sqrt{(1-1/m)} U_s)^2 \\
&= 2\left(1-\sum^{m}_{i=1} Leb(T_i \cap S_i) \right)
\end{align*}

\noindent for all $s,t \in [0,1]$.  That is, while 
$X_t$ and $\sqrt{(1-1/m)} U_t$ have different covariance 
structures, their $L^2$-structures are
identical.  The Sudakov-Fernique Inequality then 
allows us to conclude again that $\bbe 	\tilde{H}_m = \bbe D_m$ in
a manner independent of the development of Theorem \ref{thm4}.
\end{Rem}

\begin{Rem}\label{sumrmt} Let us briefly summarize the connections between
random matrix theory and
the Brownian functionals encountered in this paper.
Writing $x^{(m)} = (x_1,x_2,\dots,x_m)$, and defining
$\Delta(x) = \Pi_{1\le i < j \le m}(x_i - x_j)$
to be the Vandermonde determinant, we have the following 
facts. \\

\noindent (i) $D_m \stackrel{d}{=} \lambda_1^{(m)}$,
where $\lambda_1^{(m)}$ is the largest eigenvalue of the $m \times m$ GUE,
with the scaling taken to be such that the diagonal elements $X_{i,i}$
satisfy $\bbe X_{i,i}^2 = 1$, and the off-diagonal elements
$X_{i,j}$, for $i \ne j$, satisfy $\bbe |X_{i,j}|^2 = 1$.
Using standard random matrix results (see, {\it e.g.}, \cite{Me}),
the distribution of $D_m$, 
for all $m \ge 1$ and all $s \in \bbr$, is given by

$$\bbp(D_m \le s) = c_m \int_{A_s} e^{-\sum_{i=1}^{m}x_i^2/2} \Delta(x)^2dx^{(m)},$$

\noindent where 
$$A_s = \{x \in \bbr^m: \max_{1\le i \le m} x_i \le s\},$$

\noindent where

$$c_m^{-1} = \int_{\bbr^m} e^{-\sum_{i=1}^{m}x_i^2/2} \Delta(x)^2dx^{(m)}.$$

\noindent (ii) $\tilde{H}_m \stackrel{d}{=} \lambda_1^{(m,0)}$,
where $\lambda_1^{(m,0)}$ is the 
largest eigenvalue of the $m \times m$ 
{\it traceless} GUE, with the scaling as in (i).
Using the joint 
density of the eigenvalues of the traceless
$m \times m$ GUE \cite{Me,TW},
the distribution function 
of $\tilde{H}_m$ can also be computed directly,  
for all $m \ge 2$ and all $s \ge 0$, as

$$\bbp(\tilde{H}_m \le s) = c_m^0 \int_{A_s^0} e^{-(m/2)\sum_{i=1}^{m}x_i^2} \Delta(x)^2dx^{(m,0)},$$

\noindent where $dx^{(m,0)}$ is Lebesgue measure over
the set $\{\sum_{i=1}^m x_i = 0\}$, and where
$$A_s^0 = \{x \in \bbr^m: \max_{1\le i \le m} x_i \le s\} \cap
\left\{\sum_{i=1}^m x_i = 0\right\},$$

\noindent where

$$(c_m^0)^{-1} = \int_{\{\sum_{i=1}^m x_i = 0\}} e^{-(m/2)\sum_{i=1}^{m}x_i^2} 
\Delta(x)^2dx^{(m,0)}.$$

\noindent Note that $\tilde{H}_m$ is a.s. non-negative, and so
$\bbp(\tilde{H}_m \le s) = 0$, for all $s < 0$.\\

\noindent (iii) $J_m := \sqrt{p_{max}}\{\sqrt{1 - kp_{max}}Z_k + \tilde{H}_k\}$,
the limiting functional of Proposition~\ref{prop1} for
the $m$-letter non-uniform case,
having its most probable letters of multiplicity $k$ 
occuring with probability $p_{max}$,
is equal in law to the sum of a normal
random variable and a variable whose distribution,
up to the scaling factor $\sqrt{p_{max}}$,
is that of the largest eigenvalue of the $k \times k$ 
{\it traceless} GUE, with the scaling as in (i) and (ii).
(Note also that since
$D_m \stackrel{d}{=} Z_m + \tilde{H}_m$,
$J_m$ is also equal in law to the sum
of a normal random variable and a variable
whose distribution,
up to the scaling factor $\sqrt{p_{max}}$, 
is that of the largest 
eigenvalue of the $k \times k$ GUE.)  Its,
Tracy, and Widom \cite{ITW1} show that,
for all $m \ge 2$ and all $s \in \bbr$, 
$J_m$ has distribution given by 

$$\bbp(J_m \le s) = c_{k,p_{max}} \int_{A_s} 
   e^{-\frac{1}{2p_{max}}\left[\sum_{i=1}^{k}x_i^2 
      + \frac{p_{max}}{1-kp_{max}}(\sum_{i=1}^{k}x_i)^2\right]} \Delta(x)^2dx^{(k)},$$

\noindent where 
$$A_s = \{x \in \bbr^k: \max_{1\le i \le k} x_i \le s\},$$

\noindent and where

$$c_{k,p_{max}}^{-1} = \int_{\bbr^k} 
   e^{-\frac{1}{2p_{max}}\left[\sum_{i=1}^{k}x_i^2 
      + \frac{p_{max}}{1-kp_{max}}(\sum_{i=1}^{k}x_i)^2\right]} 
\Delta(x)^2dx^{(k)}.$$

\noindent  Moreover, in the discussion prior to 
Theorem \ref{thm2}, we noted that the
$k$-fold integral representation of the limiting
distribution of $J_m$ came from simplifying
a more complex expression.  This expression
described the distribution
of $J_m$ as that of largest eigenvalue of the
direct sum of $d$ mutually independent GUEs, 
each of size $k_j \times k_j$, $1 \le j \le d$,
subject to the eigenvalue constraint
$\sum_{i=1}^m \sqrt{p_i}\lambda_i = 0$.
Here the $k_j$ were the multiplicities of
the probabilities having common values,
the $p_i$ were ordered in decreasing order,
and the eigenvalues were ordered in terms
of the GUEs corresponding to the 
appropriate values of $p_i$.

Note that when $k=1$, the 
limiting distribution becomes simply

$$\bbp(J_m \le s) = 
    \frac{1}{ \sqrt{ 2\pi p_{max} (1-p_{max}) } } \int_{-\infty}^s 
    e^{- x^2/2p_{max}(1-p_{max}) } dx,$$

\noindent which is simply a $N(0,p_{max}(1-p_{max}))$ distribution.\\

\noindent (iv) The Tracy-Widom distribution function $F_2$, 
which also describes the limiting distribution of 
$(L\sigma_n - 2\sqrt{n})/n^{1/6}$, (see \cite{BDJ}),
is given, for all $t \in \bbr$, by

$$F_2(t) = \exp\left(-\int_t^{\infty} (x-t)u^2(x) dx \right),$$

\noindent where $u(x)$ is the solution to the Painlev\'e II
equation $u_{xx} = 2u^3 + xu$ with $u(x) \sim -Ai(x)$,
as $x \rightarrow \infty$, where $Ai$ is the Airy function.
\end{Rem}

\section{Countable Infinite Alphabets}

Let us now study the problem of describing
$LI_n$ for an ordered, countably infinite alphabet
${\cal A} =\{\alpha_n\}_{n\ge1}$, where
$\alpha_1 < \alpha_2 < \cdots < \alpha_m < \cdots$.
Let $(X_i)^{n}_{i=1}$, $X_i \in{\cal A}$,
be an iid sequence, with 
$\bbp(X_1 = \alpha_r) = p_r > 0$,
for $r \ge 1$. 
 
The central idea in the first 
part of our approach is to
introduce two new sequences derived from
$(X_i)^{n}_{i=1}$. Fix $m \ge 1$.
The first sequence, which we shall term the 
{\it capped sequence}, is defined by taking  
$T_i^m = X_i \wedge \alpha_m$, for $i \ge 1$. 
The second one, 
$(Y^m_i)^{N_{n,m}}_{i=1}$,
the {\it reduced sequence},
consists of the subsequence of
$(X_i)^{n}_{i=1}$
of length $N_{n,m}$,
for which $X_i \le \alpha_m$, 
for $i \ge 1$.
Thus, the capped sequence 
$(T^m_i)^{n}_{i=1}$ 
is obtained by setting to $\alpha_m$ all letter 
values greater than $\alpha_m$, 
while the reduced sequence
$(Y^m_i)^{N_{n,m}}_{i=1}$
is obtained by eliminating letter 
values greater than $\alpha_m$ altogether.

Let $LI^{cap}_{n,m}$ and $LI^{red}_{n,m}$ 
to be the lengths of the longest increasing subsequence
of $(T^m_i)^n_{i=1}$ and $(Y^m_i)^{N_{n,m}}_{i=1}$, 
respectively.  Now on the one hand, any subsequence 
of the reduced sequence is again a subsequence 
of the original sequence 
$(X_i)^{n}_{i=1}$.  On the other hand, any
increasing subsequence of $(X_i)^{n}_{i=1}$
is again an increasing subsequence of the capped
one.  These two observations lead to the pathwise bounds

\begin{equation}\label{ia1}
LI^{red}_{n,m} \le LI_n \le LI^{cap}_{n,m},
\end{equation}

\noindent for all $m\ge 1$ and $n \ge 1$.

These bounds suggest that the 
behavior of the iid infinite case perhaps mirrors 
that of the iid finite-alphabet case.  Indeed,
we do have the following result, which amounts
to an extension of Theorem \ref{thm2} (or, 
more precisely, of Proposition \ref{prop1}) 
to the iid infinite-alphabet case.

\begin{theorem}\label{thm6}
Let $(X_i)_{i\ge1}$ be a sequence of 
iid random variables taking values in 
the ordered alphabet
${\cal A} =\{\alpha_n\}_{n\ge1}$.
Let $\bbp(X_1 = \alpha_j) = p_j$, for
$j \ge 1$.
Let $p_{max} = p_{j_1} =p_{j_2} = \cdots \ = p_{j_k}$, 
$1 \le j_1 < j_2 < \cdots < j_k$, $k \ge 1$, and 
let $p_i < p_{max}$, otherwise. Then

\begin{align*}
\frac{LI_n-p_{max}n}{\sqrt n} \Rightarrow 
\sqrt{p_{max}}\{\sqrt{1 - p_{max}k}Z_k + H_k\}:=R(p_{max},k).  
\end{align*}
\end{theorem}

The proof of the theorem relies on an understanding
of the limiting distributions of $LI^{red}_{n,m}$ and
$LI^{cap}_{n,m}$.  To this end, let us introduce 
some more notation.  For a finite
$m$-alphabet, and for $V_1, \dots, V_n$ iid 
with $\bbp(V_1=\alpha_r) = q_r>0$, let
$LI_n(q):=LI_n(q_1,\dots,q_m)$ denote 
the length of the longest increasing
subsequence of $(V_i)^{n}_{i=1}$.  For each $m \ge 1$, let 
also $\pi_m = \sum_{r=1}^{m}p_r$.  

First, let us choose $m$
large enough so that $1 - \pi_{m-1} < p_{max}$.
Next, observe that, from the capping at $\alpha_m$,
$LI^{cap}_{n,m}$ is distributed as $LI_n(\tilde{p})$,
where $\tilde{p}=(p_1,\dots,p_{m-1},1-\pi_{m-1})$.
But since $m$ is chosen large enough, 
the maximal probability among the entries
of $\tilde{p}$ is then $p_{max}$, of multiplicity $k$,
as for the original infinite alphabet.  
By Theorem \ref{thm2}, we thus have

\begin{equation}\label{ia2}
\frac{LI_n(\tilde{p})-p_{max}n}{\sqrt n} \Rightarrow R(p_{max},k),
\end{equation}

\noindent as $n \rightarrow \infty$.

Turning to $LI^{red}_{n,m}$, suppose that the number
of elements $N_{n,m}$ of the reduced subsequence
$(Y^m_i)^{N_{n,m}}_{i=1}$ is equal to $j$.
Since only the elements of 
$(X_i)^{n}_{i=1}$ which are at most $\alpha_m$ are left, 
$LI^{red}_{n,m}$ must be distributed as $LI_j(\hat{p})$,
where $\hat{p}=(p_1/\pi_m,\dots,p_m/\pi_m)$.
From the way $m$ is chosen,
the maximal probability among the entries
of $\hat{p}$ is then $p_{max}/\pi_m$, of multiplicity $k$.
Invoking again the finite-alphabet 
result of Theorem \ref{thm2}, we find that

\begin{equation}\label{ia3}
\frac{LI_n(\hat{p})-(p_{max}/\pi_m)n}{\sqrt n} 
   \Rightarrow R\left(\frac{p_{max}}{\pi_m},k \right),
\end{equation}

\noindent as $n \rightarrow \infty$.

We now relate the two limiting expressions in
\eqref{ia2} and \eqref{ia3} by the following
elementary lemma.

\begin{lemma}\label{lem1}
Let $k \ge 1$ be an integer, and
let $(q_m)^{\infty}_{m=1}$ be a 
sequence of reals
in $(0,1/k]$ converging to $q \ge 0$. Then
$R(q_m,k) \Rightarrow R(q,k)$,
as $m \rightarrow \infty$.
\end{lemma}

\noindent \begin{Proof} 
Assume $q > 0$. Then

\begin{align}\label{ia4}
R(q_m,k) 
  &= \sqrt{q_m}\left\{\sqrt{1 - q_mk}Z_k 
     + H_k\right\}\nonumber\\
  &= \sqrt{q_m} \{ \sqrt{1 - qk}Z_k + H_k \}\nonumber\\
  &\qquad + \sqrt{q_m}
    \{ \sqrt{1 - q_mk} - \sqrt{1 - qk} \}Z_k\nonumber\\
  &= \sqrt{\frac{q_m}{q}}  R(q,k) + c_mZ_k,
\end{align}

\noindent where 
$c_m = \sqrt{1 - q_mk} - \sqrt{1 - qk}$.  Since 
$q_m \rightarrow q$ as $m \rightarrow \infty$,
$c_m \rightarrow 0$,
and so $c_mZ_k \Rightarrow 0$, as $m \rightarrow \infty$.
This gives the result. 
The degenerate case, $q = 0$, is clear.\CQFD

\end{Proof}

The main idea developed in the proof of
Theorem \ref{thm6} is now to use the basic inequality
\eqref{ia1} in conjunction 
with a conditioning argument
for $LI^{red}_{n,m}$, in order
to apply Lemma \ref{lem1},
{\it i.e.,} to use
$R(p_{max}/\pi_m,k) \Rightarrow R(p_{max},k),$
as $m \rightarrow \infty$,
since $\pi_m \rightarrow 1$.\\

\noindent \begin{Proof} {\bf (Theorem \ref{thm6})} 
First, fix an arbitrary $s > 0$. 
As previously noted in Remark \ref{rem1},
$R(p_{max},k)$ has a density supported on $\bbr$
($\bbr^{+}$ in the uniform case), and so $s$ is a 
continuity point
of its distribution function.   
Next, choose $0 < \epsilon_1 < 1$, and $0 < \delta < 1$,
and again note that $(1+\delta)s$ is also necessarily
a continuity point for $R(p_{max},k)$.

With this choice of $\epsilon_1$, pick $\beta > 0$
such that $\bbp(Z \ge \beta) < \epsilon_1/2$,
where $Z$ is a standard normal random variable.
Finally, pick $\epsilon_2$ such that
$0 < \epsilon_2 < \epsilon_1 \bbp(R(p_{max},k) < (1+\delta)s)$.
Such a choice of $\epsilon_2$ can always be made
since the support of $R(p_{max},k)$ includes $\bbr^{+}$.

We have seen that, for $m$ large enough,
we can bring some finite-alphabet results to bear
on the infinite case.  In fact, we need a few more
technical requirements on $m$ to complete our proof.
Setting $\sigma^2_m = \pi_m(1-\pi_m)$,
we choose large enough $m$ so that:

\begin{align*}
\text{(i)}   &\quad1 - \pi_{m-1} < p_{max},\\
\text{(ii)}  &\quad (s + p_{max}\beta\sigma_m/\pi_m)/\sqrt{\pi_m - \beta\sigma_m} < (1+\delta)s, \text{ and}\\
\text{(iii)} &\quad\left| \bbp(R(p_{max},k) < (1+\delta)s)
                    - \bbp(R(p_{max}/\pi_m,k) < (1+\delta)s)\right| < \epsilon_2/2.
\end{align*}

\noindent The conditions (i) and (ii) 
are clearly satisfied, since
$\pi_m \rightarrow 1$  and
$\sigma_m \rightarrow 0$, as $m \rightarrow \infty$.
The condition (iii) is also satisfied, 
as seen by applying Lemma \ref{lem1}
to $R(p_{max}/\pi_m,k)$,
with $\pi_m \rightarrow 1$,
and since $(1+\delta)s$ is also
a continuity point for $R(p_{max},k)$.

Now recall that  
$LI^{cap}_{n,m}$ is distributed as $LI_n(\tilde{p})$,
where $\tilde{p}=(p_1,\dots,p_{m-1},1-\pi_{m-1})$.
Hence, we have from \eqref{ia1} 
and \eqref{ia2} that

\begin{align}\label{ia5}
\frac{LI_n - p_{max}n}{\sqrt{n}} &\le \frac{LI^{cap}_{n,m}- p_{max}n} {\sqrt{n}}
  \Rightarrow R(p_{max},k),
\end{align}

and so

\begin{align}\label{ia5b}
\bbp\left(\frac{LI_n - p_{max}n}{\sqrt{n}} \le s\right) 
 &\ge \bbp\left(\frac{LI^{cap}_{n,m}- p_{max}n} {\sqrt{n}} \le s\right)\nonumber\\
 &\rightarrow \bbp(R(p_{max},k) \le s),
\end{align}

\noindent as $n \rightarrow \infty$ (and, in fact, for any $s\in \bbr$).

More work is required to 
make use of the left-most minorization
in \eqref{ia1} ({\it i.e.}, $LI^{red}_{n,m} \le LI_n$.)
Recall that if the length 
$N_{n,m}$ of the reduced sequence is
equal to $j$, then $LI^{red}_{n,m}$ must be 
distributed as $LI_j(\hat{p})$,
where $\hat{p}=(p_1/\pi_m,\dots,p_m/\pi_m)$.
Now the essential observation is that $N_{n,m}$
is distributed as a binomial random variable with
parameters $\pi_m$ and $n$.  It is thus natural 
to focus on the values of $j$ 
close to $\bbe N_{n,m} = n\pi_m$.
Writing the variance of $N_{n,m}$ as
$n\sigma^2_m$, where, as above, $\sigma^2_m = \pi_m(1-\pi_m)$,
and 
$$\gamma_{n,m,j} := \bbp(N_{n,m} = j ) = \left( \begin{array}{c}
n  \\
j
\end{array} \right) \pi_m^j(1-\pi_m)^{n-j},$$
\noindent we have

{\allowdisplaybreaks
\begin{align}
&\bbp\left( \frac{LI^{red}_{n,m}- p_{max}n} {\sqrt{n}} \le s \right) \nonumber\\
&= \sum_{j=0}^n \bbp\left( \frac{LI^{red}_{n,m}- p_{max}n} {\sqrt{n}} \le s \vert
       N_{n,m}=j \right) \gamma_{n,m,j}\nonumber\\
&= \sum_{j=0}^n \bbp\left( \frac{LI_j(\hat{p})-p_{max}n}{\sqrt n} \le s \right) 
       \gamma_{n,m,j}\nonumber\\
&= \sum_{j=0}^n \bbp\left( \frac{LI_j(\hat{p})-\frac{p_{max}}{\pi_m}j}{\sqrt j} 
       \le \sqrt{\frac{n}{j}}\left( s +\frac{p_{max}}{\sqrt{n}}
       \left(n-\frac{j}{\pi_m}\right) \right)  \right)  
       \gamma_{n,m,j}\nonumber\\
&\le \sum_{j=\lceil n\pi_m - \beta\sigma_m\sqrt{n}\rceil}^{n}
       \bbp\left( \frac{LI_j(\hat{p})-\frac{p_{max}}{\pi_m}j}{\sqrt j} 
       \le \sqrt{\frac{n}{j}}\left( s +\frac{p_{max}}{\sqrt{n}}
       \left(n-\frac{j}{\pi_m}\right) \right)  \right)  
       \gamma_{n,m,j}\nonumber\\
& \qquad \qquad + \sum_{j=0}^{\lceil n\pi_m - \beta\sigma_m\sqrt{n}\rceil-1}
       \gamma_{n,m,j}\nonumber\\
&< \sum_{j=\lceil n\pi_m - \beta\sigma_m\sqrt{n}\rceil}^{n}
       \bbp\left( \frac{LI_j(\hat{p})-\frac{p_{max}}{\pi_m}j}{\sqrt j} 
       \le \sqrt{\frac{n}{j}}\left( s +\frac{p_{max}}{\sqrt{n}}
       \left(n-\frac{j}{\pi_m}\right) \right) \right)  
       \gamma_{n,m,j}\nonumber\\
& \qquad \qquad + \epsilon_1,\label{ia6}
\end{align}}

\noindent for sufficiently large $n$,
where \eqref{ia6} follows from
the Central Limit Theorem and our choice of $\beta$,
and where, as usual, $\lceil\cdot\rceil$ 
is the ceiling function.

Next, note that for 
$\lceil n\pi_m - \beta\sigma_m\sqrt{n}\rceil \le j \le n$,
and by condition (ii),

\begin{align}\label{ia7}
&\sqrt{\frac{n}{j}}\left( s +\frac{p_{max}}{\sqrt{n}}
       \left(n-\frac{j}{\pi_m}\right) \right) \nonumber\\
&\qquad\qquad <   \sqrt{\frac{n}{n\pi_m - \beta\sigma_m\sqrt{n}}}\left( s +\frac{p_{max}}{\sqrt{n}}
       \left(n-\frac{n\pi_m - \beta\sigma_m\sqrt{n}}{\pi_m}\right) \right) \nonumber\\
&\qquad\qquad =   \frac{1}{\sqrt{\pi_m - \beta\sigma_m/\sqrt{n}}}
        \left( s +\frac{p_{max}\beta\sigma_m}{\pi_m} \right) \nonumber\\
&\qquad\qquad \le \frac{1}{\sqrt{\pi_m - \beta\sigma_m}}
        \left( s +\frac{p_{max}\beta\sigma_m}{\pi_m} \right) \nonumber\\
&\qquad\qquad < s(1+\delta).
\end{align}

Hence, for sufficiently large $n$, we have

\begin{align}\label{ia8}
&\sum_{j=\lceil n\pi_m - \beta\sigma_m\sqrt{n}\rceil}^{n}
       \bbp\left( \frac{LI_j(\hat{p})-\frac{p_{max}}{\pi_m}j}{\sqrt j} 
       \le \sqrt{\frac{n}{j}}\left( s +\frac{p_{max}}{\sqrt{n}}
       \left(n-\frac{j}{\pi_m}\right) \right) \right)  
       \gamma_{n,m,j}\nonumber\\
& \qquad \qquad + \epsilon_1\nonumber\\
& \qquad \le \sum_{j=\lceil n\pi_m - \beta\sigma_m\sqrt{n}\rceil}^{n}
       \bbp\left( \frac{LI_j(\hat{p})-\frac{p_{max}}{\pi_m}j}{\sqrt j} 
       \le s(1+\delta)\right)  
       \gamma_{n,m,j} 
       + \epsilon_1.
\end{align}

Now from the condition (iii),
and from the weak convergence,
as $j \rightarrow \infty$, of 
$(LI_j(\hat{p})-(p_{max}/\pi_m)j)/\sqrt{j}$
to $R(p_{max}/\pi_m,k)$, 
we find that, 
for $j$ large enough,

\begin{align}\label{ia9}
&\left| \bbp \left(\frac{LI_j(\hat{p})-\frac{p_{max}}{\pi_m}j}{\sqrt j} 
\le (1+\delta)s \right)
      - \bbp ( R(p_{max},k) \le (1+\delta)s  )  \right|
       \nonumber\\
&\le  \left| \bbp \left(\frac{LI_j(\hat{p})-\frac{p_{max}}{\pi_m}j}{\sqrt j} 
       \le (1+\delta)s  \right)
      - \bbp \left( R\left(\frac{p_{max}}{\pi_m},k\right) \le (1+\delta)s  \right)  \right|
      \nonumber\\
&\quad + \left|\bbp ( R(p_{max},k) \le (1+\delta)s  ) 
             - \bbp \left( R\left(\frac{p_{max}}{\pi_m},k\right) \le (1+\delta)s  
          \right)  \right|
          \nonumber\\
&< \frac{\epsilon_2}{2} + \frac{\epsilon_2}{2} \nonumber\\
&< \epsilon_1 \bbp ( R(p_{max},k) \le (1+\delta)s  ),
\end{align}

\noindent and so,

\begin{equation}\label{ia10}
\bbp \left(\frac{LI_j(\hat{p})-\frac{p_{max}}{\pi_m}j}{\sqrt j} \le (1+\delta)s \right)
\le (1+\epsilon_1)\bbp \left( R(p_{max},k) \le (1+\delta)s  \right).
\end{equation}

Now since
$\lceil n\pi_m - \beta\sigma_m\sqrt{n} \rceil \rightarrow \infty$,
as $n \rightarrow \infty$,
with the help of
\eqref{ia8} and \eqref{ia10},
\eqref{ia6} becomes

\begin{align}\label{ia11}
&\bbp\left( \frac{LI^{red}_{n,m}- p_{max}n} {\sqrt{n}} \le s \right) \nonumber\\
&\quad \quad\le \sum_{j=\lceil n\pi_m - \beta\sigma_m\sqrt{n}\rceil}^{n}
       (1+\epsilon_1)\bbp \left( R(p_{max},k) \le (1+\delta)s  \right) 
       \gamma_{n,m,j} + \epsilon_1\nonumber\\
&\quad \quad\le (1+\epsilon_1)\bbp \left( R(p_{max},k) \le (1+\delta)s  \right)
            + \epsilon_1.
\end{align}

From \eqref{ia1} we know that 
$LI^{red}_{n,m} \le LI_n$ a.s., and so

\begin{align}\label{ia12}
\bbp\left(\frac{LI_n - p_{max}n}{\sqrt{n}} \le s\right) 
 &\le \bbp\left(\frac{LI^{red}_{n,m}- p_{max}n} {\sqrt{n}} \le s\right)\nonumber\\ 
 & \le (1+\epsilon_1)\bbp \left( R(p_{max},k) \le (1+\delta)s  \right)
            + \epsilon_1,
\end{align}

\noindent for large enough $n$.
But since $\epsilon_1$ and $\delta$ are arbitrary,
\eqref{ia12} and \eqref{ia5b} together show that

\begin{align}\label{ia13}
\bbp\left( \frac{LI_n - p_{max}n}{\sqrt{n}} \le s \right) 
&\rightarrow \bbp ( R(p_{max},k) \le s),
\end{align}

\noindent for all $s > 0$.

The proof for $s < 0$ is similar.  Indeed, since
necessarily $p_{max} < 1/k$, $R(p_{max},k)$ describes 
the limiting distribution of the longest increasing
subsequence for a non-uniform alphabet,
and so is supported on $\bbr$.  One then needs only 
examine quantities of the form, {\it e.g.},
$\bbp (R(p_{max},k) \le (1-\delta)s)$, instead of
$\bbp (R(p_{max},k) \le (1+\delta)s)$, as we have done
throughout the proof for $s > 0$.  These changes lead to the 
resulting statement.\CQFD \end{Proof}

\begin{Rem}\label{rem1b}
As an alternative to the above proof, 
one could certainly adopt the 
finite-alphabet development
of the previous sections so as to express 
$LI_n$, for countable infinite alphabets, 
in terms of approximations
to functionals of Brownian motion.
More precisely,

\begin{align*}
LI_n &= \sup_{m \ge 2} \max_{\stackrel{\scriptstyle 0\le k_1\le\cdots}{\le k_{m-1}\le n}}
\left\{S^1_{k_1}+S^2_{k_2}+ \cdots + S^{m-1}_{k_{m-1}}+a^{m}_n\right\}\nonumber\\
&= \sup_{m \ge 2} \left\{\frac {n}{m}-\frac1m\sum^{m-1}_{r=1} rS^r_n+
\max_{\stackrel{\scriptstyle 0\le k_1\le\cdots}{\le k_{m-1}\le n}}
\sum^{m-1}_{r=1} S^r_{k_r}\right\},
\end{align*}

\noindent where $a^m_n$ counts the number of occurrences
of the letter $\alpha_m$ among $(X_i)_{1\le i\le n}$,
and $S^r_k=\sum^k_{i=1}Z^r_i$ is the sum of independent
random variables defined as in \eqref{item3}.  
After centering and normalizing the $S^r_k$, as was done to obtain
\eqref{item8} in the non-uniform finite alphabet development, 
one could then try to apply Donsker's Theorem to obtain a Brownian
functional, which we now know to be distributed as $R(p_{max},k$).

\end{Rem}

\section{Concluding Remarks}

Our development of the general finite-alphabet case 
leads us to consider several new directions
in which to pursue this method and raises a number 
of interesting questions.  These include the following.\\

\noindent \textbullet \quad Extending our fixed finite-alphabet 
case to that of having each $X_n$ take values in
$\{1,2,\dots,m_n\}$ is an important first step.  
Fruitful approaches to such asymptotic questions would nicely
close the circle of ideas initiated here.  Such a study is already
under consideration (see \cite{HL3}).

\noindent \textbullet \quad  As we have noted throughout the paper, 
there is a pleasing if still rather mysterious connection between
our limiting distribution results and those of random matrix theory. 
This connection deserves to be further explored. 
Recall, for instance, Baryshnikov's
observation \cite{Ba} that the process $D_m$ is identical in law to the process
$\lambda^{(m)}_1$ consisting of the largest eigenvalues of 
the $m \times m$ minor of an infinite GUE matrix. 
This fact is consistent with an interleaving-eigenvalue result from 
basic linear algebra, namely, that if $\lambda_1 \ge \lambda_2 \ge
\dots \ge \lambda_n$ are the eigenvalues of an $n \times n$ 
symmetric matrix $A$, and if $\mu_1 \ge \mu_2 \ge \dots \ge \mu_{n-1}$
are the eigenvalues of the matrix consisting of the 
first $(n-1)$ rows and columns of $A$, then $\lambda_1 \ge \mu_1 \ge \lambda_2 \ge
\dots \ge \mu_{n-1} \ge \lambda_n$.   We thus see the 
consistency between the $D_m \le D_{m+1}$ a.s. fact noted above and
that of $\lambda_1 \ge \mu_1$.  

\noindent \textbullet \quad  Pursuing our analysis further, 
one might hope to find ways in which we can 
derive the densities 
of our limiting functionals in a direct manner.  
Its, Tracy, and Widom \cite{ITW1} have obtained clear expressions 
of the limiting distributions.  While we have 
obtained our limiting distributions 
in a rather direct way, in turn, these densities 
do not clearly follow from our approach.  This point deserves more work.  

\noindent \textbullet \quad  In another direction, 
our independent-letter paradigm can be 
extended to various types of dependent
cases, foremost of which would be the Markov case.  
This will be presented elsewhere \cite{HL2},
where the framework of \cite{HLM} is, moreover,
further extended.   

\noindent \textbullet \quad  Various other types of subsequence 
problems can be tackled by the methodologies used 
in the present paper.  To name but a few, comparisons for 
unimodal sequences, 
alternating sequences, and sequences with 
blocks will deserve further similar studies.\\

{\bf Acknowledgements} C.H.~would like to 
thank Zhan Shi for discussions and encouragements on this project.
Both authors would like to thank the organizers of
the Special Program on High-Dimensional Inference
and Random Matrices at SAMSI.  
Their hospitality and support, 
through the grant DMS-0112069, 
greatly facilitated the
completion of this paper.

\noindent File reference: finalpha40.tex/pdf

\end{document}